\documentclass[12pt]{article}

\usepackage{latexsym,amsfonts,amsmath,amssymb,amsthm,url}
\usepackage[english]{babel}
\usepackage[latin1]{inputenc}
\usepackage[T1]{fontenc}
\usepackage[dvips]{graphicx}
\usepackage[dvips]{color}

\makeatletter

\@addtoreset{equation}{section}
\makeatother

\begin{document}

\newtheorem{lemm}{Lemma}[section]
\newtheorem{rema}[lemm]{Remark}
\newtheorem{prop}[lemm]{Proposition}
\newtheorem{theo}[lemm]{Theorem}
\newtheorem{fact}[lemm]{Fact}
\newtheorem{coro}[lemm]{Corollary}
\newtheorem*{ack}{Acknowledgments}

\newcommand{\p}{\ensuremath{\bar{p}}}
\newcommand{\q}{\ensuremath{\bar{q}}}
\renewcommand{\P}{\ensuremath{\mathbf{P}}}
\newcommand{\E}{\ensuremath{\mathbf{E}}}
\newcommand{\N}{\ensuremath{\mathbb{N}}}
\newcommand{\Z}{\ensuremath{\mathbb{Z}}}
\newcommand{\R}{\ensuremath{\mathbb{R}}}
\newcommand{\B}{\ensuremath{B}}
\newcommand{\Cste}[1]{\ensuremath{c_{#1}}}
\newcommand{\defeg}{\ensuremath{\overset{\hbox{\tiny{def}}}{=}}}
\newcommand{\Indic}{{\mathchoice {\rm 1\mskip-4mu l} {\rm 1\mskip-4mu l} {\rm 1\mskip-4.5mu l} {\rm 1\mskip-5mu l}}}

\title{\textbf{Rate of growth of a transient cookie random walk} }

\author{
\normalsize{\textsc{Anne-Laure Basdevant} and \textsc{Arvind Singh}
\footnote{Address for both authors: \it{Laboratoire de Probabilités
et Modèles Aléatoires, Université Pierre et Marie Curie, 175 rue du
Chevaleret, 75013 Paris, France.} }}}

\date{}

\maketitle \vspace{-1cm}

\begin{center}
University Paris VI
\end{center}

\vspace*{0.2cm}

\begin{abstract}
We consider a one-dimensional transient cookie random walk. It is
known from a previous paper \cite{BasdevantSingh06-preprint} that a
cookie random walk $(X_n)$ has positive or zero speed according to
some positive parameter $\alpha >1$ or $\le 1$. In this  article, we
give the exact rate of growth of $(X_n)$ in the zero speed regime,
namely: for $0< \alpha <1$, $X_n/n^{\frac{\alpha+1}{2}}$ converges
in law to a Mittag-Leffler distribution whereas for $\alpha=1$, $
X_n(\log n)/n$ converges in probability to some positive constant.
\end{abstract}

\bigskip
{\small{

\bigskip
\noindent{\bf A.M.S. Classification. }60K35, 60J80, 60F05

\bigskip
\noindent{\bf e-mail. } anne-laure.basdevant@ens.fr,
arvind.singh@ens.fr }}

\section{Introduction}\label{ART5_sectionintro}

Let us pick a strictly positive integer $M$. An $M$-cookie random
walk (also called multi-excited random walk) is a walk on $\Z$ which
has a bias to the right upon its $M$ first visits at a given site
and evolves like a symmetric random walk afterwards. This model was
introduced by Zerner \cite{Zerner06} as a generalization, in the
one-dimensional setting, of the model of the excited random walk
studied by Benjamini and Wilson \cite{BenjaminiWilson03}. In this paper,
we consider the case where the initial cookie environment is
spatially homogeneous.  Formally, let $(\Omega,\P)$ be some
probability space and choose a vector $\p = (p_1,\ldots, p_M)$ such
that $p_i \in [\frac{1}{2},1)$ for all $i=1,\ldots,M$. We say that
$p_i$ represents the strength of the $i^{\hbox{\tiny{th}}}$ cookie
at a given site. Then, an $(M,\p)$-cookie random walk $(X_n,\,
n\in\N)$ is a nearest neighbour random walk, starting from $0$, and
with transition probabilities:
\begin{equation*}
\P\{X_{n+1} = X_n + 1\, \hbox{|} \, X_0,\ldots,X_n\} = \left\{
\begin{array}{ll}
p_j&\hbox{if $j = \sharp\{0\leq i\leq n,\, X_i=X_n\}\leq M$,}\\
\frac{1}{2}&\hbox{otherwise.}
\end{array}
\right.
\end{equation*}
In particular, the future position $X_{n+1}$ of the walk after time
$n$ depends on the whole trajectory  $X_0,X_1,\ldots,X_n$.
Therefore, $X$ is not, unless in degenerated cases, a Markov
process. The cookie random walk is a rich stochastic model.
Depending on the cookie environment $(M,\p)$, the process can either
be transient or recurrent. Precisely, Zerner \cite{Zerner06} (who
considered an even more general setting) proved, in our case, that
if we define
\begin{equation}\label{ART5_defalpha}
\alpha = \alpha(M,\p)  \defeg \sum_{i=1}^{M}(2p_i - 1) - 1,
\end{equation}
then
\begin{itemize}
\item if $\alpha\leq 0$, the cookie random walk is recurrent,
\item if $\alpha>0$, the cookie random walk is transient towards $+\infty$.
\end{itemize}
Thus, a $1$-cookie random walk is always recurrent but, for two or
more cookies, the walk can either be transient or recurrent. Zerner
also proved that the limiting velocity of the walk is well defined.
That is, there exists a deterministic constant $v=v(M,\p)\geq 0$
such that
\begin{equation*}
\lim_{n\to\infty}\frac{X_n}{n} = v\quad\hbox{almost surely.}
\end{equation*}
However, we may have $v=0$. Indeed, when there are at most two
cookies per site, Zerner proved that $v$ is always zero. On the
other hand, Mountford \emph{et al.} \cite{MountfordPimentelValle06}
showed that it is possible to have $v>0$ if the number of cookies is
large enough. In a previous paper \cite{BasdevantSingh06-preprint},
the authors showed that, in fact, the strict positivity of the speed
depends on the position of $\alpha$ with respect to $1$:
\begin{itemize}
\item if $\alpha\leq 1$, then $v=0$,
\item if $\alpha> 1$, then $v>0$.
\end{itemize}
In particular, a positive speed may be obtained with just three
cookies per site. The aim of this paper is to find the exact
rate of growth of a transient cookie random walk in zero speed regime.
In this perspective, numerical simulations of
Antal and Redner \cite{AntalRedner05} indicated that, for a transient
2-cookies random walk, the expectation of $X_n$ is of order $n^\nu$,
for some constant $\nu\in (\frac{1}{2},1)$ depending on the strength
of the cookies. We shall prove that, more generally, $\nu =
\frac{\alpha+1}{2}$.

\begin{theo}\label{ART5_MainTheo}
Let $X$ be a $(M,\p)$-cookie random walk and let $\alpha$ be defined
by \textup{(\ref{ART5_defalpha})}. Then, when the walk is transient
with zero speed, \emph{i.e.} when $0<\alpha\leq 1$,
\begin{enumerate}
\item If $\alpha < 1$,
\begin{equation*}
\frac{X_n}{n^{\frac{\alpha+1}{2}}}\underset{n\to\infty}
{\overset{\hbox{\tiny{law}}}{\longrightarrow}} \mathcal{M}_{\frac{\alpha+1}{2}}
\end{equation*}
where $\mathcal{M}_{\frac{\alpha+1}{2}}$ denotes a Mittag-Leffler
distribution with parameter $\frac{\alpha+1}{2}$.
\item If $\alpha=1$, there exists a constant $c>0$ such that
\begin{equation*}
\frac{\log
n}{n}X_n\underset{n\to\infty}{\overset{\hbox{\tiny{prob.}}}{\longrightarrow}}
c.
\end{equation*}
\end{enumerate}
These results also hold with $\sup_{i\leq n}X_i$ and $\inf_{i\geq
n}X_i$ in place of $X_n$.
\end{theo}
\newpage

\begin{figure}[here]
\setlength\unitlength{1cm}
\begin{picture}(15,7)
\put(2,0){\includegraphics[height=7cm,angle=0]{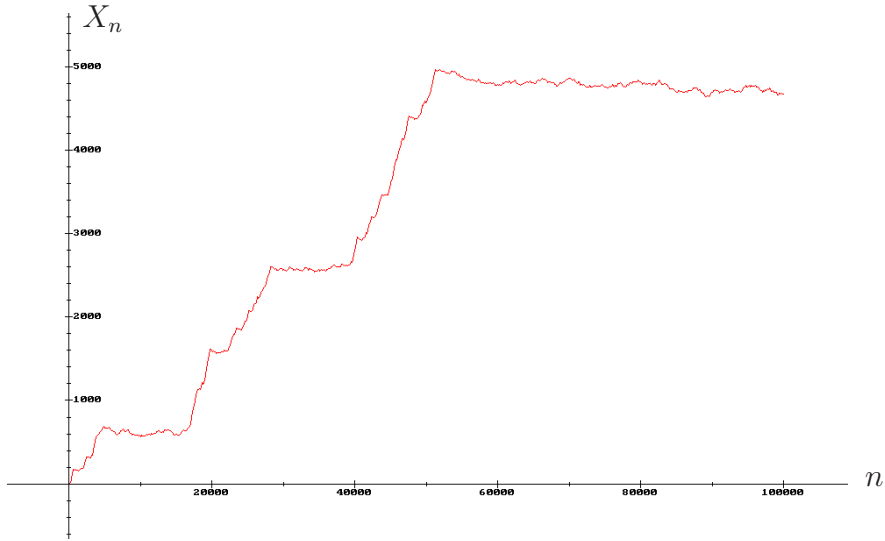}}
\put(13.4,0.7){$n$} \put(3,6.8){$X_n$}
\end{picture}
\caption{Simulation of the $100000$ first steps of a cookie random
walk with $M=3$ and $p_1=p_2=p_3=\frac{3}{4}$ (\emph{i.e.} $\alpha
=\frac{1}{2}$ and $\nu=\frac{3}{4}$).}
\end{figure}

This theorem bears many likenesses to the famous result of
Kesten-Kozlov-Spitzer \cite{KestenKozlovSpitzer75} concerning the
rate of transience of a one-dimensional random walk in random
environment. Indeed, following the method initiated in
\cite{BasdevantSingh06-preprint}, we can reduce the study of the
walk to that of an auxiliary Markov process $Z$. In our setting, $Z$
is a branching process with migration. By comparison, Kesten
\emph{et al.} obtained the rates of transience of the random walk in
random environment via the study of an associated branching process
in random environment. However, the process $Z$ considered here and
the process introduced in \cite{KestenKozlovSpitzer75} have quite
dissimilar behaviours and the methods used for their study are
fairly different.

Let us also note that, as $\alpha$ tends to zero, the rate of growth
$n^{(1+\alpha)/2}$ tends to $\sqrt{n}$. This suggests that, when the
cookie walk is recurrent (\emph{i.e.} $-1<\alpha\leq 0$), its growth
should not be much larger than that of a simple symmetric random
walk. In fact, we believe that, in the recurrent setting,
$\sup_{i\leq n}X_i$ should be of order $l(n)\sqrt{n}$ for some
slowly varying function $l$.

The remainder of this paper is organized as follow. In the next
section, we recall the construction of the associated process $Z$
described in \cite{BasdevantSingh06-preprint} as well as some
important results concerning this process. In section
\ref{ART5_sectionsigma}, we study the tail distribution of the
return time to zero of the process $Z$. Section
\ref{ART5_sectionprogeny} is devoted to estimating the tail
distribution of the total progeny of the branching process over an
excursion away from $0$. The proof of this result is based on
technical estimates whose proofs are given in section
\ref{ART5_sectiontechnical}. Once all these results obtained, the
proof of the main theorem is quite straightforward and is finally
given in the last section.

\section{The process $Z$}\label{ART5_section2}

In the rest of this paper, $X$ will denote an $(M,\p)$-cookie random
walk. We will also always assume that we are in the transient regime
and that the speed of the walk is zero, that is
\begin{equation*}
0<\alpha\leq 1.
\end{equation*}
The proof of Theorem \ref{ART5_MainTheo} is based on a careful study of
the hitting times of the walk:
\begin{equation*}
T_n \defeg \inf\{k\geq 0,\, X_k = n\}.
\end{equation*}
We now introduce a Markov process $Z$ closely connected with these
hitting times. Indeed, we can summarize Proposition 2.2 and equation
(4) of \cite{BasdevantSingh06-preprint} as follows:
\begin{prop}\label{ART5_propZT} There exist a Markov process $(Z_n,\, n\in\N)$ starting from $0$ and a
sequence of random variables $(K_n,\, n\geq 0)$ converging in law
towards a finite random variable $K$ such that, for each $n$
\begin{equation*}
T_n \overset{\hbox{\tiny{law}}}{=} n + 2\sum_{k=0}^{n}Z_k + K_n.
\end{equation*}
\end{prop}
Therefore, a careful study of $Z$ will enable us to obtain precise
estimates on the distribution of the hitting times. In the rest of
this section, we shall recall the construction of $Z$ and some
important results obtained in \cite{BasdevantSingh06-preprint}.

For each $i=1,2,\ldots$, let $B_i$ be a Bernoulli random variable
with distribution
\begin{equation*}
\P\{B_i = 1\} = 1 - \P\{B_i = 0\} = \left\{
\begin{array}{ll}
p_i&\hbox{ if $1\leq i\leq M$,}\\
\frac{1}{2}&\hbox{ if $i > M$.}
\end{array}\right.
\end{equation*}
We define the random variables $A_0,A_1,\ldots,A_{M-1}$ by
\begin{equation*}
A_j \defeg \sharp\{1\leq i \leq k_j,\, B_i = 0\} \quad
\hbox{where}\quad k_j \defeg \inf\Big(i\geq 1,\, \sum_{l=1}^{i} B_l =
j+1\Big).
\end{equation*}
Therefore, $A_j$ represents the number of "failures" before having
$j+1$ "successes"  along the sequence of coin tossings $(B_i)$. It
is to be noted that the random variables $A_j$ admit some
exponential moments:
\begin{equation}\label{ART5_expmomA}
\E[s^{A_j}]<\infty\quad\hbox{ for all $s\in[0,2)$.}
\end{equation}
According to Lemma $3.3$ of \cite{BasdevantSingh06-preprint}, we also have
\begin{equation}\label{ART5_espA}
\E[A_{M-1}] = 2\sum_{i=1}^{M}(1-p_i) = M - 1 - \alpha.
\end{equation}
Let $(\xi_i,\, i\in\N^*)$ be a sequence of i.i.d. geometric random
variables with parameter $\frac{1}{2}$ (\emph{i.e.} with mean $1$),
independent of the $A_j$. The process $Z$ mentioned above is a
Markov process with transition probabilities given by
\begin{equation}\label{ART5_defProcessZ}
\P\big\{ Z_{n+1} = j\,\hbox{|}\,Z_n = i\big\} = \P\Big\{
\Indic_{\{i\leq M-1\}}A_{i} + \Indic_{\{i > M-1\}}\Big(A_{M-1} +
\sum_{k=1}^{i- M + 1}\xi_k\Big) = j\Big\}.
\end{equation}
As usual, we will use the notation $\P_x$ to describe the law of the
process starting from $x\in\N$ and $\E_x$ the associated
expectation, with the conventions $\P = \P_0$ and $\E = \E_0$. Let us
notice that $Z$ may be interpreted as a branching process with
random migration, that is, a branching process which allows both
immigration and emigration components.
\begin{itemize}
\item If $Z_n = i \in \{M,M+1,\ldots\}$, then  $Z_{n+1}$ has the law of
$\sum_{k=1}^{i-M+1}\xi_k + A_{M-1}$, \emph{i.e.} $M-1$ particles
emigrate from the system and the remaining particles reproduce
according to a geometrical law with parameter $\frac{1}{2}$ and
there is also an immigration of $A_{M-1}$ new particles.
\item If $Z_n = i \in \{0,\ldots,M-1\}$, then $Z_{n+1}$ has the same law as
$A_i$, \emph{i.e.} all the $i$ particles emigrate the system and
$A_i$ new particles immigrate.
\end{itemize}
Since we assume that the cookie vector $\p$ is such that $p_i<1$ for
all $i$, the process $Z$ is an irreducible Markov process. More
precisely,
\begin{equation*}
\P_x\{Z_1 = y\} >0\quad\hbox{for all $x,y\in\N$.}
\end{equation*}
From the construction of the random variables $A_i$, we have $A_0
\leq A_1 \leq \ldots \leq A_{M-1}$. This fact easily implies that,
for any $x\leq y$, the process $Z$ under $\P_x$ (starting from $x$)
is stochastically dominated by $Z$ under $\P_y$ (starting from $y$).
Let us also note that, for any $k\geq M-1$,
\begin{equation}\label{ART5_martineg}
 \E[Z_{n+1}-Z_{n}\,|\,Z_n = k] =  \E[A_{M-1}] -M + 1 = -\alpha.
\end{equation}
This quantity is negative and we say that emigration dominates
immigration. In view of (\ref{ART5_martineg}), a simple martingale
argument shows that $Z$ is recurrent. More precisely, according to
section $2$ of \cite{BasdevantSingh06-preprint}, the process $Z$ is,
in fact, positive recurrent and thus converges in law, independently
of its starting point, towards a random variable $Z_\infty$ whose
law is the unique invariant probability for $Z$. Moreover, according
to Remark $3.7$ of \cite{BasdevantSingh06-preprint}, the tail
distribution of $Z_\infty$ is regularly varying with index $\alpha$:
\begin{prop}\label{ART5_proptailZinfty} There exists a constant $c>0$ such that
\begin{equation*}
\P\{Z_\infty> x\}\underset{x\to\infty}{\sim} \left\{
\begin{array}{ll}
c/x^\alpha&\hbox{if $\alpha\in(0,1)$,}\\
c\log x /x &\hbox{if $\alpha = 1$.}
\end{array}\right.
\end{equation*}
\end{prop}
Let now $\sigma$ denote the first return time to $0$ for the process
$Z$,
\begin{equation*}
\sigma \defeg \inf\{n\geq 1,\, Z_n = 0\}.
\end{equation*}
According to the classical expression of the invariant probability,
for any non negative function $f$, we have
\begin{equation}\label{ART5_eqexc}
\E\!\left[\sum_{i=0}^{\sigma-1}f(Z_i)\right] =
\E[\sigma]\E[f(Z_\infty)].
\end{equation}
In particular, we deduce the following corollary which will be found very
useful:
\begin{coro}\label{ART5_CorIntegZinfty} We have, for $\beta\geq 0$,
\begin{equation*}
\E\!\left[\sum_{i=0}^{\sigma-1}Z^\beta_i\right] \left\{
\begin{array}{ll}
< \infty &\hbox{if $\beta < \alpha$,}\\
= \infty &\hbox{if $\beta \geq \alpha$.}\\
\end{array}\right.
\end{equation*}
\end{coro}
\section{The return time to zero}\label{ART5_sectionsigma}
We have already stated that $Z$ is an irreducible positive recurrent
Markov chain, thus the return time $\sigma$ to zero has finite
expectation. The aim of this section is to strengthen this result by
giving the asymptotic of the tail distribution of $\sigma$.
Precisely, we will show that
\begin{prop}\label{ART5_PropTailSigma}
For any initial starting point $x\geq 1$, there exists $c = c(x) >0$
such that
\begin{equation*}
\P_{x}\{\sigma > n\}
\underset{n\to\infty}{\sim}\frac{c}{n^{\alpha+1}}.
\end{equation*}
\end{prop}
Notice that we do not allow the starting point $x$ to be $0$. In
fact, this assumption could be dropped but it would unnecessarily
complicate the proof of the proposition which is technical enough
already. Yet, we have already mentioned that $Z$ starting from $0$
is stochastically dominated by $Z$ starting from $1$, thus
$\P\{\sigma>n\}\leq \P_1\{\sigma>n\}$. We also have $\P\{\sigma>n\}
\geq \P\{Z_1=1\}\P_1\{\sigma>n-1\}$. Therefore, we deduce that
\begin{equation*}
\frac{c_1}{n^{\alpha+1}}\leq \P\{\sigma>n\} \leq
\frac{c_2}{n^{\alpha+1}}
\end{equation*}
where $c_1$ and $c_2$ are two strictly positive constants. In
particular, we obtain the following corollary which will be
sufficient for our needs.
\begin{coro}\label{ART5_approxsigma}
We have
\begin{equation}
\E[\sigma^{\beta}]  \left\{
\begin{array}{ll}
< \infty&\hbox{ if $\beta < \alpha+1$,}\\
= \infty&\hbox{ if $\beta \geq \alpha+1$.}
\end{array}
\right.
\end{equation}
\end{coro}
The method used in the proof of the proposition is classical and
based on the study of probability generating functions. Proposition
\ref{ART5_PropTailSigma} was first proved by Vatutin
\cite{Vatutin77} who considered a branching process with exactly one
emigrant at each generation. This result was later generalized for
branching processes with more than one emigrant by Vinokurov
\cite{Vinokurov87} and also by Kaverin \cite{Kaverin90}. However, in
our setting, we deal with a branching process with migration, that
is, where both immigration and emigration are allowed. More
recently, Yanev and Yanev proved similar results for such a class of
processes, under the assumption that, either there is at most one
emigrant per generation \cite{YanevYanev04} or that immigration
dominates emigration \cite{YanevYanev95} (in our setting, this would
correspond to the case $\alpha<0$).

For the process $Z$, the emigration component dominates the
immigration component and this leads to some additional technical
difficulties. Although there is a vast literature on the subject
(see the authoritative survey of Vatutin and Zubkov
\cite{VatutinZubkov93} for additional references), we did not find a
proof of Proposition \ref{ART5_PropTailSigma} in our setting. We
shall therefore provide here a complete argument but we invite the
reader to look in the references mentioned above for additional
details.

Recall the definition of the random variables $A_i$ and $\xi_i$ defined
in section \ref{ART5_section2}. We introduce, for $s\in[0,1]$,
\begin{eqnarray*}
F(s) &\defeg& \E[s^{\xi_{1}}] = \frac{1}{2-s},\\
\delta(s) &\defeg& (2-s)^{M-1}\E[s^{A_{M-1}}], \\
H_k(s) &\defeg&
(2-s)^{M-1-k}\E[s^{A_{M-1}}]-\E[s^{A_k}]\quad\hbox{for $1\leq k\leq
M-2$.}
\end{eqnarray*}
Let $F_j(s) \defeg F\circ \ldots\circ F(s)$ stand for the $j$-fold of
$F$ (with the convention $F_0 = \hbox{Id}$). We also define by
induction
\begin{equation*}
\left\{
\begin{array}{l}
\gamma_0(s) \defeg 1,\\
\gamma_{n+1}(s) \defeg \delta(F_n(s))\gamma_n(s).
\end{array}\right.
\end{equation*}
We use the abbreviated notations $F_j \defeg F_j(0)$, $\gamma_n
\defeg \gamma_n(0)$. We start with a simple lemma.
\begin{lemm}\label{ART5_lemmeS1}
\begin{enumerate}
\item[\textup{(a)}] $F_n = 1 - \frac{1}{n+1}$.
\item[\textup{(b)}] $H_k(1-s) = -H'_k(1) s + \mathcal{O}(s^2)$ when $s\to 0$ for all $1\leq k\leq M-2$.
\item[\textup{(c)}] $\delta(1-s) = 1+\alpha s + \mathcal{O}(s^2)$ when $s\to 0$.
\item[\textup{(d)}] $\gamma_n \sim_{\infty} c_3 n^\alpha$ with $c_3>0$.
\end{enumerate}
\end{lemm}
\begin{proof}
Assertion (a) is straightforward. According to (\ref{ART5_expmomA}), the
functions $H_k$ are analytic on $(0,2)$ and (b) follows from a
Taylor expansion near $1$. Similarly, (c) follows from a Taylor
expansion near $1$ of the function $\delta$ combined with
(\ref{ART5_espA}). Finally, $\gamma_n$ can be expressed in the form
\begin{equation*}
\gamma_n = \prod_{j=0}^{n-1}\delta(F_j)
\underset{n\to\infty}{\sim}\prod_{j=1}^{n}\left(1 +
\frac{\alpha}{j}\right)\underset{n\to\infty}{\sim} c_3 n^\alpha,
\end{equation*}
which yields (d).
\end{proof}
Let $\widetilde{Z}$ stand for the process $Z$ absorbed at $0$:
\begin{equation*}
\widetilde{Z}_n \defeg Z_{n}\Indic_{\{n\leq\inf(k,\, Z_k=0)\}}.
\end{equation*}
We also define, for $x\geq 1$ and $s\in[0,1]$,
\begin{eqnarray}
\label{ART5_refGG}J_x(s) &\defeg& \sum_{i=0}^{\infty}\P_{x}\{\widetilde{Z}_i\neq 0\}s^i,\\
\nonumber G_{n,x}(s) &\defeg& \E_x[s^{\widetilde{Z}_n}],
\end{eqnarray}
and for $1\leq k\leq M-2$,
\begin{equation*}
g_{k,x}(s) \defeg \sum_{i=0}^{\infty}\P_{x}\{\widetilde{Z}_i =
k\}s^{i+1}.
\end{equation*}
\begin{lemm}\label{ART5_lemmeS2} For any $1\leq k\leq M-2$, we have
\begin{enumerate}
\item[\textup{(a)}] $\sup_{x\geq 1} g_{k,x}(1) < \infty$.
\item[\textup{(b)}] for all $x\geq 1$, $g'_{k,x}(1) <\infty.$
\end{enumerate}
\end{lemm}
\begin{proof}
The value $g_{x,k}(1)$ represents the expected number of visits to
site $k$ before hitting $0$ for the process $Z$ starting from $x$.
Thus, an easy application of the Markov property yields

\begin{equation*}
g_{k,x}(1) = \frac{\P_x\{\hbox{$Z$ visits $k$ before
$0$}\}}{\P_k\{\hbox{$Z$ visits $0$ before returning to $k$\}}} <
\frac{1}{\P_{k}\{Z_1 = 0\} } < \infty.
\end{equation*}
This proves (a). We now introduce the return times $\sigma_k
\defeg \inf(n\geq 1,\, Z_n =k)$. In view of the Markov property, we
have
\begin{eqnarray*}
g'_{k,x}(1) &=& g_{k,x}(1) + \E_{x}\Big[\sum_{n=1}^{\infty}n\Indic_{\{\widetilde{Z}_n=k\}}\Big]\\
&=& g_{k,x}(1) + \sum_{i=1}^{\infty}\P_{x}\{\sigma_k = i,\,
\sigma_k <\sigma\}\E_{k}\Big[\sum_{n=0}^{\infty}(i+n)\Indic_{\{\widetilde{Z}_n=k\}}\Big]\\
&=& g_{k,x}(1) + \E_{x}[\sigma_k\Indic_{\{\sigma_k <\sigma\}
}]g_{k,k}(1) + \P_{x}\{\sigma_k
<\sigma\}\E_{k}\Big[\sum_{n=0}^{\infty}n\Indic_{\{\widetilde{Z}_n=k\}}\Big].
\end{eqnarray*}
Since $Z$ is a positive recurrent Markov process, we have
$\E_{x}[\sigma_k\Indic_{\{\sigma_k <\sigma\} }] \leq \E_{x}[\sigma]
<\infty$. Thus, it simply remains to show that
$\E_{k}\Big[\sum_{n=0}^{\infty}n\Indic_{\{\widetilde{Z}_n=k\}}\Big]<\infty$.
Using the Markov property, as above, but considering now the partial
sums, we get, for any $N\geq 1$,
\begin{eqnarray*}
\E_{k}\left[\sum_{n=1}^{N}n\Indic_{\{\widetilde{Z}_n=k\}}\right] &=&
\sum_{i=1}^{N}\P_{k}\{\sigma_k =i,\, \sigma_k < \sigma\}
\E_{k}\left[\sum_{n=0}^{N-i}(i+n)\Indic_{\{\widetilde{Z}_n=k\}}\right]\\
&\leq& \E_k\left[\sigma_k\Indic_{\{\sigma_k<\sigma\}}\right]g_{k,k}(1)
+\P_k\{\sigma_k<\sigma\}\E_{k}\left[\sum_{n=1}^{N}n\Indic_{\{\widetilde{Z}_n=k\}}\right].
\end{eqnarray*}
Since $\P_k\{\sigma < \sigma_k\} \geq \P_{k}\{Z_1 = 0\} >0$, we
deduce that
\begin{equation*}
\E_{k}\left[\sum_{n=1}^{N}n\Indic_{\{\widetilde{Z}_n=k\}}\right] \leq
\frac{\E_k\left[\sigma_k\Indic_{\{\sigma_k<\sigma\}}\right]g_{k,k}(1)}{\P_k\{\sigma < \sigma_k\}}<\infty.
\end{equation*}
and we conclude the proof letting $N$ tend to $+\infty$.
\end{proof}

\begin{lemm}\label{ART5_lemmeS3} The function $J_x$ defined by (\ref{ART5_refGG})
may be expressed in the form
\begin{equation*}
J_x(s) =  \widehat{J}_x(s) +
\sum_{k=1}^{M-2}\widetilde{J}_{k,x}(s)\quad\hbox{for $s\in[0,1)$},
\end{equation*}
where
\begin{equation*}
\widehat{J}_x(s) \defeg \frac{\sum_{n=0}^{\infty}\gamma_n (1
-(F_n)^x) s^n}{(1-s)\sum_{n=0}^\infty \gamma_n s^n}\quad\hbox{ and
}\quad \widetilde{J}_{k,x}(s) \defeg
\frac{g_{k,x}(s)\sum_{n=0}^{\infty}\gamma_n H_k(F_n)
s^n}{(1-s)\sum_{n=0}^\infty \gamma_n s^n}.
\end{equation*}
\end{lemm}
\begin{proof} From the definition (\ref{ART5_defProcessZ}) of the branching
process $Z$, we get, for $n\geq 0$,
\begin{equation*}
\begin{aligned}
&G_{n+1,x}(s) = \E_{x}\Big[\E_{\widetilde{Z}_n}[s^{\widetilde{Z}_1}]\Big]\\
&\!\!= \P_x\{\widetilde{Z}_n=0\} +
\sum_{k=1}^{M-2}\P_x\{\widetilde{Z}_n = k\}\E[s^{A_k}]
+ \sum_{k=M-1}^{\infty}\P_x\{\widetilde{Z}_n = k\}\E[s^{\xi}]^{k-(M-1)}\E[s^{A_{M-1}}]\\
&\!\!=\!\left(\!\!1\!-\!\frac{\E[s^{A_{M-1}}]}{\E[s^{\xi}]^{M-1}}\right)\!\P_x\{\widetilde{Z}_n=0\}
\!-\! \sum_{k=1}^{M-2}\!\P_x\{\widetilde{Z}_n = k\}H_k(s) +
\frac{\E[s^{A_{M-1}}]}{\E[s^{\xi}]^{M-1}}\sum_{k=0}^{\infty}\P_x\{\widetilde{Z}_n
= k\}\E[s^{\xi}]^{k}\!.
\end{aligned}
\end{equation*}
Since $\E[s^{\xi}] = F(s)$ and $G_{n,x}(0) = \P_x\{\widetilde{Z}_n=
0\}$, using the notation introduced in the beginning of the section,
the last equality may be rewritten
\begin{equation*}
G_{n+1,x}(s) = \delta(s)G_{n,x}(F(s)) + (1-\delta(s))G_{n,x}(0) -
\sum_{k=1}^{M-2}\P_x\{\widetilde{Z}_n = k\}H_k(s).
\end{equation*}
Iterating this equation then setting $s=0$ and using the relation
$G_{0,x}(F_{n+1}) = (F_{n+1})^x$, we deduce that, for any $n\geq 0$,
\begin{equation}\label{ART5_step2}
G_{n+1,x}(0) = \sum_{i=0}^{n}(1\!-\!\delta(F_i))\gamma_i
G_{n-i,x}(0) \,+\, \gamma_{n+1}(F_{n+1})^x -
\sum_{k=1}^{M-2}\sum_{i=0}^{n}\P_{x}\{\widetilde{Z}_{n-i}=k\}\gamma_i
H_k(F_{i}).
\end{equation}
Notice also that $\P_x\{\widetilde{Z}_n \neq 0\} = 1 - G_{n,x}(0)$.
In view of (\ref{ART5_step2}) and making use of the relation
$(1-\delta(F_i))\gamma_i = \gamma_i - \gamma_{i+1}$, we find, for
all $n\geq 0$ (with the convention $\sum_0^{-1}=0$)
\begin{eqnarray*}
\P_x\{\widetilde{Z}_n \neq 0\} &=& \gamma_n(1-(F_n)^x) + \sum_{i=0}^{n-1}(\gamma_i -
\gamma_{i+1})\P_x\{\widetilde{Z}_{n-1-i} \neq 0\}\\
&& +
\sum_{k=1}^{M-2}\sum_{i=0}^{n-1}\P_{x}\{\widetilde{Z}_{n-1-i}=k\}\gamma_i
H_k(F_{i}).
\end{eqnarray*}
Therefore, summing over $n$, for $s<1$,
\begin{equation*}
\begin{aligned}
&J_x(s) = \sum_{n=0}^{\infty}\P_x\{\widetilde{Z}_n \neq 0\} s^n\\
&= \sum_{n=0}^{\infty}\gamma_n(1-(F_n)^x)s^n + \sum_{n=0}^{\infty}
\sum_{i=0}^{n}(\gamma_i - \gamma_{i+1})\P_x\{\widetilde{Z}_{n-i} \neq 0\}s^{n+1}\\
&
\qquad + \sum_{k=1}^{M-2} \sum_{n=0}^{\infty}
\sum_{i=0}^{n}\P_{x}\{\widetilde{Z}_{n-i}=k\}\gamma_i
H_k(F_{i})s^{n+1}\\
&= \sum_{n=0}^{\infty}\gamma_n (1-(F_n)^x)s^n +
J_x(s)\sum_{n=0}^{\infty}(\gamma_n - \gamma_{n+1})s^{n+1} +
\sum_{k=1}^{M-2} g_{k,x}(s)\sum_{n=0}^{\infty}\gamma_nH_k(F_n)s^{n}.
\end{aligned}
\end{equation*}
We conclude the proof noticing that $\sum_{n=0}^{\infty}(\gamma_n -
\gamma_{n+1})s^{n+1} = (s-1)\sum_{n=0}^{\infty}\gamma_n s^n + 1$.
\end{proof}

We can now give the proof of the proposition.

\begin{proof}[Proof of Proposition \ref{ART5_PropTailSigma}]
Recall that the parameter $\alpha$ is such that $0<\alpha\leq 1$. We
first assume $\alpha<1$. Fix $x\geq 1$ and $1\leq k\leq M-2$. In
view of Lemma \ref{ART5_lemmeS1} and with the help of an
Abelian/Tauberian theorem (\emph{c.f.} Chap VIII of
\cite{Feller71}), we check that
\begin{equation*}
(1-s)\sum_{n=0}^\infty \gamma_n s^n \underset{s\to 1^-}{\sim}
\frac{c_3\Gamma(\alpha+1)}{(1-s)^\alpha}\quad\hbox{ and }\quad
\sum_{n=0}^{\infty}\gamma_n H_k(F_n) s^n \underset{s\to 1^-}{\sim}
-\frac{c_3 H'_k(1)\Gamma(\alpha)}{(1-s)^\alpha}.
\end{equation*}
These two equivalences show that $\widetilde{J}_{k,x}(1)\defeg
\lim_{s\to 1^-}\widetilde{J}_{k,x}(s)$ is finite. More precisely, we
get
\begin{equation*}
\widetilde{J}_{k,x}(1) =  -\frac{g_{k,x}(1)H'_k(1)}{\alpha},
\end{equation*}
so that we may write
\begin{equation}\label{ART5_sty}
\frac{\widetilde{J}_{k,x}(1) - \widetilde{J}_{k,x}(s)}{1-s} = \left(\frac{g_{k,x}(1) -
g_{k,x}(s)}{1-s}\right)\frac{\widetilde{J}_{k,x}(s)}{g_{k,x}(s)} + \frac{g_{k,x}(1)
\widetilde{B}_k(s)}{(1-s)^2\sum_{n=0}^\infty \gamma_n s^n}
\end{equation}
with the notation
\begin{equation*}
\widetilde{B}_k(s) \defeg \frac{H'_k(1)}{\alpha}(s-1)\sum_{n=0}^\infty \gamma_n s^n -
\sum_{n=0}^{\infty}\gamma_n H_k(F_n) s^n.
\end{equation*}
The first term on the r.h.s. of (\ref{ART5_sty}) converges towards
$-g'_k(1)H'_k(1)/\alpha$ as $s$ tends to $1$ (this quantity is
finite thanks to Lemma \ref{ART5_lemmeS2}). Making use of the relation
$\gamma_{n+1} = \delta(F_n)\gamma_n$, we can also rewrite
$\widetilde{B}_k$ in the form
\begin{equation*}
\widetilde{B}_k(s) =
\sum_{n=1}^{\infty}\gamma_{n-1}\left[\frac{H'_k(1)}{\alpha}(1-\delta(F_{n-1}))
- \delta(F_{n-1})H_k(F_n)\right]s^n - \frac{H'_k(1)}{\alpha}-H_k(0).
\end{equation*}
With the help of Lemma \ref{ART5_lemmeS1}, it is easily check that
\begin{equation*}
\gamma_{n-1}\left[\frac{H'_k(1)}{\alpha}(1-\delta(F_{n-1})) -
\delta(F_{n-1})H_k(F_n)\right] =
\mathcal{O}\left(\frac{1}{n^{2-\alpha}}\right).
\end{equation*}
Since  $\alpha<1$, we conclude that
\begin{equation}\label{ART5_refyop2}
\widetilde{B}_k(1) = \lim_{s\to 1^-}\widetilde{B}_k(s) \quad\hbox{is finite.}
\end{equation}
We also have
\begin{equation}\label{ART5_refyop3}
(1-s)^2\sum_{n=0}^\infty \gamma_n s^n \underset{s\to 1^-}{\sim}
\frac{c_3\Gamma(\alpha+1)}{(1-s)^{\alpha-1}}.
\end{equation}
Thus, combining (\ref{ART5_sty}), (\ref{ART5_refyop2}) and (\ref{ART5_refyop3}), as
$s\to 1^-$,
\begin{equation}\label{ART5_refJWT}
\frac{\widetilde{J}_{k,x}(1) - \widetilde{J}_{k,x}(s)}{1-s} =
\frac{g_{k,x}(1)
\widetilde{B}_k(1)}{c_3\Gamma(\alpha+1)}(1-s)^{\alpha-1} +
o\big((1-s)^{\alpha-1}\big).
\end{equation}
We can deal with $\widehat{J}_x$ in exactly the same way. We now
find $\widehat{J}_x(1) = \frac{x}{\alpha}$ and setting
\begin{equation}\label{ART5_refBBX}
\widehat{B}_x(1) \defeg
\sum_{n=1}^{\infty}\gamma_{n-1}\left[\frac{x}{\alpha}(\delta(F_{n-1})-1)
- \delta(F_{n-1})(1-(F_n)^x)\right] + \frac{x}{\alpha}- 1,
\end{equation}
we also find that, as $s\to 1^-$,
\begin{equation}\label{ART5_refJWH}
\frac{\widehat{J}_{x}(1) - \widehat{J}_{x}(s)}{1-s} =
\frac{\widehat{B}_x(1)}{c_3\Gamma(\alpha+1)}(1-s)^{\alpha-1} +
o\big((1-s)^{\alpha-1}\big).
\end{equation}
Putting together (\ref{ART5_refJWT}) and (\ref{ART5_refJWH}) and using Lemma
\ref{ART5_lemmeS3}, we obtain
\begin{equation}\label{ART5_mlf1}
\frac{J_{x}(1) - J_{x}(s)}{1-s} = C_x (1-s)^{\alpha-1} +
o\big((1-s)^{\alpha-1}\big)
\end{equation}
with
\begin{equation}\label{ART5_refCC}
C_x \defeg \frac{1}{c_3\Gamma(\alpha+1)}\left(\widehat{B}_x(1) +
\sum_{k=1}^{M-2}g_{k,x}(1)\widetilde{B}_k(1)\right).
\end{equation}
Since $x\neq 0$, we have $\P_{x}\{\widetilde{Z}_n \neq 0\} =
\P_{x}\{\sigma >n\}$ and, from the definition of $J_x$, we deduce
\begin{equation}\label{ART5_ml2}
\sum_{n=0}^{\infty}\Big(\sum_{k=n+1}^{\infty}\P_{x}\{\sigma
>k\}\Big)s^n = \frac{J_x(1)-J_x(s)}{1-s}.
\end{equation}
Combining (\ref{ART5_mlf1}) and (\ref{ART5_ml2}), we see that $C_x \geq 0$.
Moreover, the use of two successive Tauberian theorems yields
\begin{equation*}
\P_{x}\{\sigma > n\} = \frac{C_x \alpha}{\Gamma(1-\alpha)
n^{\alpha+1}} + o\!\left(\frac{1}{n^{\alpha+1}}\right).
\end{equation*}
It remains to prove that $C_x \neq 0$. To this end, we first notice
that, for $x,y\geq 0$, we have $\P_y\{Z_1 = x\}>0$ and
\begin{equation*}
\P_y\{\sigma > n\} \geq \P_y\{Z_1 = x\}\P_{x}\{\sigma > n-1\}.
\end{equation*}
Thus, $C_y\geq \P_y\{Z_1 = x\} C_x$ so it suffices to show that
$C_x$ is not zero for some $x$. In view of (a) of Lemma
\ref{ART5_lemmeS2}, the quantity
\begin{equation*}
\sum_{k=1}^{M-2}g_{k,x}(1)\widetilde{B}_k(1)
\end{equation*}
is bounded in $x$. Looking at the expression of $C_x$ given in
(\ref{ART5_refCC}),  it just remains to prove that $\widehat{B}_x(1)$ can
be arbitrarily large. In view of (\ref{ART5_refBBX}), we can write
\begin{equation*}
\widehat{B}_x(1) = x S(x) + \frac{x}{\alpha}- 1
\end{equation*}
where
\begin{equation*}
S(x) \defeg
\sum_{n=1}^{\infty}\gamma_{n-1}\left[\frac{1}{\alpha}(\delta(F_{n-1})-1)
- \delta(F_{n-1})\frac{(1-(F_n)^x)}{x}\right].
\end{equation*}
But for each fixed $n$, the function
\begin{equation*}
x\to \delta(F_{n-1})\frac{(1-(F_n)^x)}{x}
\end{equation*}
 decreases to $0$ as $x$ tends to infinity, so the monotone
convergence theorem yields
\begin{equation*}
S(x) \underset{x\to\infty}{\uparrow}
\sum_{n=1}^{\infty}\frac{\gamma_{n-1}}{\alpha}(\delta(F_{n-1})-1)
\;\sim\; c_3\sum_{n=1}^{\infty}\frac{1}{n^{1-\alpha}} \;=\; +\infty.
\end{equation*}
Thus, $\widehat{B}_x(1)$ tends to infinity as $x$ tends to infinity
and the proof of the proposition for $\alpha<1$ is complete. The
case $\alpha=1$ may be treated in a similar fashion (and it is even
easier to prove that the constant is not zero). We skip the details.
\end{proof}

\begin{rema}\label{ART5_remarksigma}
The study of the tail distribution of the return time is the key to
obtaining conditional limit theorems for the branching process, see
for instance \cite{Kaverin90,Vatutin77,Vinokurov87,YanevYanev04}.
Indeed, following Vatutin's scheme \cite{Vatutin77} and using
Proposition \ref{ART5_PropTailSigma}, it can now be proved that
$Z_n/n$ conditioned on not hitting 0 before time $n$ converges in
law towards an exponential distribution. Precisely, for each $x=
1,2,\ldots$ and $r\in\R_+$,
\begin{equation*}
\lim_{n\to\infty}\P_{x}\left\{\frac{Z_n}{n} \leq r\hbox{\textup{ |
}}\sigma >n \right\} = 1 - e^{-r}.
\end{equation*}
It is to be noted that this result is exactly the same as that
obtained for a classical critical Galton-Watson process (\emph{i.e.}
when there is no migration). Although, in our setting, the return
time to zero has a finite expectation, which is not the case for the
critical Galton-Watson process, the behaviours of both processes
conditionally on their non-extinction are still quite similar.
\end{rema}
\section{Total progeny over an excursion}\label{ART5_sectionprogeny}
The aim of this section is to study the distribution of the total
progeny of the branching process $Z$ over an excursion away from
$0$. We will constantly use the notation
\begin{equation*}
\nu \defeg \frac{\alpha+1}{2}.
\end{equation*}
In particular, $\nu$ ranges through $(\frac{1}{2},1]$. The main
result of this section is the key to the proof of Theorem
\ref{ART5_MainTheo} and states as follows.
\begin{prop}\label{ART5_progeny}There exists a constant $c>0$ such that
\begin{equation*}
 \P\left\{\sum_{k=0}^{\sigma-1} Z_k>x\right\}\underset{x\rightarrow \infty}{\sim}\left\{
 \begin{array}{ll}
c/x^{\nu}& \mbox{ if }\;\alpha\in(0,1)\\
 c\log x/x& \mbox{ if
}\; \alpha=1 .
\end{array}\right.
\end{equation*}
\end{prop}

Let us first give an informal explanation for this polynomial decay
with exponent $\nu$. In view of Remark \ref{ART5_remarksigma}, we can
expect the shape of a large excursion away from zero of the process
$Z$ to be quite similar to that of a Galton-Watson process. Indeed,
if $H$ denotes the height of an excursion of $Z$ (and $\sigma$
denotes the length of the excursion), numerical simulations show
that, just as in the case of a classical branching process without
migration, $H\approx \sigma$ and the total progeny
$\sum_{k=0}^{\sigma-1} Z_k$ is of the same order as $H\sigma$. Since
the decay of the tail distribution of $\sigma$ is polynomial with
exponent $\alpha+1$, the tail distribution of $\sum_{k=0}^{\sigma-1}
Z_k$ should then decrease with exponent $\frac{\alpha+1}{2}$.
In a way, this proposition tells us that the shape of an excursion is
very "squared".

Although there is a vast literature on the subject of branching
processes, it seems that there has not been much attention given to
the total progeny of the process. Moreover, the classical machinery
of generating functions and analytic methods, often used as a rule
in the study of branching processes seems, in our setting,
inadequate for the study of the total progeny.

The proof of Proposition \ref{ART5_progeny} uses a somewhat
different approach and is mainly based on a martingale argument. The
idea of the proof is fairly simple but, unfortunately, since we are
dealing with a discrete time model, a lot of additional technical
difficulties appear and the complete argument is quite lengthy. For
the sake of clarity, we shall first provide the skeleton of the
proof of the proposition, while postponing the proof of the
technical estimates to section \ref{ART5_sectionpreuvemun}.

Let us also note that, although we shall only study the particular
branching process associated with the cookie random walk, the method
presented here could be used to deal with a more general class of
branching processes with migration.

We start with an easy lemma stating that $\P\{\sum_{k=0}^{\sigma-1}
Z_k > x\}$ cannot decrease much faster than $\frac{1}{x^{\nu}}$.
\begin{lemm}\label{ART5_ordreprogeny}
For any $\beta>\nu$, we have
\begin{equation*}\E\left[\Big(\sum_{k=0}^{\sigma-1} Z_k\Big)^\beta\right] =
\infty.
\end{equation*}
\end{lemm}
\begin{proof} When $\alpha=\nu=1$, the result is a direct consequence of Corollary \ref{ART5_CorIntegZinfty}
of section \ref{ART5_section2}. We now assume $\alpha<1$. Hölder's
inequality gives
\begin{equation*}
\sum_{n=0}^{\sigma-1} Z_n^\alpha \le \sigma^{1-\alpha} (
\sum_{n=0}^{\sigma-1} Z_n )^\alpha.
\end{equation*}
Taking the expectation and applying again Hölder's inequality, we
obtain, for $\varepsilon >0$ small enough
\begin{equation*}
\E\left[\sum_{n=0}^{\sigma-1} Z^\alpha_n\right]\leq \E[\sigma^{ 1+
\alpha - \varepsilon}]^{\!\frac{1}{p}}\,\E\left[(
\sum_{n=0}^{\sigma-1} Z_n )^{\alpha q}\right]^{\frac{1}{q}},
\end{equation*}
with $p=\frac{1+\alpha-\varepsilon}{1-\alpha}$ and $\alpha
q=\frac{1+\alpha-\varepsilon}{2-\varepsilon/\alpha}$. Moreover,
Corollary \ref{ART5_CorIntegZinfty} states that
$\E[\sum_{n=0}^{\sigma-1} Z_n^\alpha]= \infty$ and thanks to
Corollary \ref{ART5_approxsigma}, $\E[\sigma^{ 1+ \alpha -
\varepsilon}] < \infty$. Therefore,
\begin{equation*}
\E\left[(\sum_{n=0}^{\sigma-1} Z_n)^{\alpha q}\right] =
\E\left[(\sum_{n=0}^{\sigma-1} Z_n)^{\nu +
\varepsilon'}\right]=\infty.
\end{equation*}
This result is valid for any $\varepsilon'$ small enough and
completes the proof of the lemma.
\end{proof}

\begin{proof}[Proof of Proposition \ref{ART5_progeny}]
Let us first note that, in view of an Abelian/Tauberian theorem,
Proposition \ref{ART5_progeny} is equivalent to
\begin{equation*}
\E\left[1-e^{-\lambda \sum_{k=0}^{\sigma-1} Z_k}\right]
\underset{\lambda\rightarrow 0^+}{\sim}\left\{
\begin{array}{ll}
C\lambda^{\nu}&\hbox{ if $\alpha\in(0,1)$,}\\
C\lambda \log \lambda&\hbox{ if $\alpha=1$,}\\
\end{array}\right.
\end{equation*}
where $C$ is a positive constant. We now  construct a martingale in
the following way. Let $K_\nu$ denote the modified Bessel function
of second kind with parameter $\nu$. For $\lambda>0$, we define
\begin{equation}\label{ART5_defphilanbda1}
\phi_\lambda(x)\defeg(\sqrt{\lambda}x)^{\nu}K_\nu(\sqrt{\lambda}x),\quad
\hbox{for $x>0$}.
\end{equation}
We shall give some important properties of $\phi_\lambda$ in section
\ref{ART5_sectionchoixphi}. For the time being, we simply recall that
$\phi_\lambda$ is an analytic, positive, decreasing function on
$(0,\infty)$ such that $\phi_\lambda$ and $\phi'_\lambda$ are
continuous at 0 with
\begin{equation}\label{ART5_valphi0}
\phi_\lambda(0)=2^{\nu-1}\Gamma(\nu)\quad \mbox{ and }\quad \phi'_\lambda(0)=0.
\end{equation}
Our main interest in $\phi_\lambda$ is that it satisfies the
following differential equation, for $x>0$:
\begin{equation}\label{ART5_equadiff}
-\lambda
x\phi_{\lambda}(x)-\alpha\phi'_\lambda(x)+x\phi''_{\lambda}(x)=0.
\end{equation}
Now let $(\mathcal{F}_{n},n\ge 0)$ denote the natural filtration of
the branching process $Z$ \emph{i.e.}
$\mathcal{F}_{n}\defeg\sigma(Z_k,0\le k\le n)$ and define, for
$n\geq 0$ and $\lambda>0$,
\begin{equation}\label{ART5_defWn}
W_n\defeg\phi_\lambda(Z_n)e^{-\lambda \sum_{k=0}^{n-1} Z_k }.
\end{equation}
Setting
\begin{equation}\label{ART5_defmu}
\mu(n)\defeg\E[W_{n}-W_{n+1}\;|\; \mathcal{F}_n],
\end{equation}
it is clear that the process
\begin{equation*}
Y_n\defeg W_n+\sum_{k=0}^{n-1} \mu(k)
\end{equation*}
is an $\mathcal{F}$-martingale. Furthermore, this martingale has
bounded increments since
\begin{equation*}
|Y_{n+1}-Y_n|\le |W_{n+1}-W_n|+|\mu(n)|\le
4||\phi_{\lambda}||_{\infty}.
\end{equation*}
Therefore, the use of the optional sampling theorem is legitimate
with any stopping time  with finite mean. In particular, applying
the optional sampling theorem with the first return time to $0$, we
get
\begin{equation*}\phi_{\lambda}(0)\E[e^{-\lambda
\sum_{k=0}^{\sigma-1} Z_k}] =\phi_{\lambda}(0)-
\E[\sum_{k=0}^{\sigma-1} \mu(k)],
\end{equation*} which we may be rewritten, using that $\phi_\lambda(0) = 2^{\nu-1}\Gamma(\nu)$,
\begin{equation}\label{ART5_alienZmu}\E[1-e^{-\lambda
\sum_{k=0}^{\sigma-1}
Z_k}]=\frac{1}{2^{\nu-1}\Gamma(\nu)}\E[\sum_{k=0}^{\sigma-1}
\mu(k)].\end{equation}
The proof of Proposition \ref{ART5_progeny} now
relies on a careful study of the expectation of
$\sum_{k=0}^{\sigma-1} \mu(k)$. To this end, we shall decompose
$\mu$ into several terms using a Taylor expansion of $\phi_\lambda$.
We first need the following lemma:

\newpage
\begin{lemm}\label{ART5_espzn+1-zn} $\hbox{ }$
\begin{enumerate}
\item[\textup{(a)}] There exists a function $f_1$ with $f_1(x)=0$ for all $x\ge M-1$ such that
\begin{equation*}
\E[Z_{n+1}-Z_n\;|\; \mathcal{F}_n]=-\alpha+f_1(Z_n).
\end{equation*}
\item[\textup{(b)}] There exists a function $f_2$ with $f_2(x)=f_2(M-1)$ for all $x\ge M-1$ such that
\begin{equation*}
\E[(Z_{n+1}-Z_n)^2\;|\; \mathcal{F}_n]=2Z_n+2f_2(Z_n).
\end{equation*}
\item[\textup{(c)}] For $p\in \N^*$, there exists a constant $D_p$ such that
\begin{equation*}
\E[|Z_{n+1}-Z_n|^p\;|\; \mathcal{F}_n]\le D_p
(Z^{p/2}_n+\Indic_{\{Z_n=0\}}).
\end{equation*}
\end{enumerate}
\end{lemm}
\begin{proof} Assertion (a) is just a rewriting of equation (\ref{ART5_martineg}).
Recall the notations introduced in section \ref{ART5_section2}. Recall in
particular that $\E[A_{M-1}]=M-1-\alpha$. Thus, for $j\ge M-1$, we
have
\begin{eqnarray*} \E[(Z_{n+1}-Z_n)^2\;|\;
Z_n=j]&=& \E\big[\big(A_{M-1}+\xi_1+\ldots+\xi_{j-M+1}-j\big)^2\big]\\
&=&\E\Big[\Big(\alpha+(A_{M-1}-\E[A_{M-1}])+\sum_{k=1}^{j-M+1}(\xi_k-\E[\xi_k])\Big)^2\Big]\\
&=&
\alpha^2+\mathbf{V}\hbox{ar}(A_{M-1})+(j-M+1)\mathbf{V}\hbox{ar}(\xi_1)\\
&=& 2Z_n+\alpha^2+\mathbf{V}\hbox{ar}(A_{M-1})-2(M-1).
\end{eqnarray*}
This proves (b). When $p$ is an even integer, we have
$\E[|Z_{n+1}-Z_n|^p\;|\; \mathcal{F}_n]=\E[(Z_{n+1}-Z_n)^p\;|\;
\mathcal{F}_n]$ and assertion (c) can be proved by developing
$(Z_{n+1}-Z_n)^p$ in the same manner as for (b). Finally, when $p$
is an odd integer, Hölder's inequality gives
\begin{equation*}
\E[|Z_{n+1}-Z_n|^p\;|\; Z_n=j>0]\le \E[|Z_{n+1}-Z_n|^{p+1}\;|\;
Z_n=j>0]^{\frac{p}{p+1}}\le
D_{p+1}^{\frac{p}{p+1}}Z_n^{\frac{p}{2}}.
\end{equation*}
\end{proof}

\noindent\emph{Continuation of the proof of Proposition
\ref{ART5_progeny}}. For $n\in [1,\sigma-2]$, the random variables $Z_n$
and $Z_{n+1}$ are both non zero and, since $\phi_\lambda$ is
infinitely differentiable on $(0,\infty)$, a Taylor expansion yields
\begin{equation}\label{ART5_taylorformul}
\phi_\lambda(Z_{n+1})
=\phi_\lambda(Z_{n})+\phi'_\lambda(Z_{n})(Z_{n+1}-Z_n) +
\frac{1}{2}\phi''_\lambda(Z_{n})(Z_{n+1}-Z_n)^2+\theta_n,
\end{equation}
where $\theta_n$ is given by Taylor's integral remainder formula
\begin{equation}\label{ART5_deftheta}
\theta_n
\defeg (Z_{n+1}-Z_n)^2\int_0^1(1-t)(\phi_\lambda''(Z_n+t(Z_{n+1}-Z_n))-\phi_\lambda''(Z_n))dt.
\end{equation}
When $n=\sigma-1$, this result is \emph{a priori} incorrect because
then $Z_{n+1}=0$. However, according to (\ref{ART5_valphi0}) and
(\ref{ART5_equadiff}), the functions $\phi_\lambda(t)$,
$\phi_\lambda'(t)$ and $t\phi_\lambda''(t)$ have finite limits as
$t$ tends to $0^+$, thus equation (\ref{ART5_taylorformul}) still holds
when $n = \sigma-1$. Therefore, for $n\in[1,\sigma-1]$,
\begin{multline*}\E[e^{\lambda
Z_n}\phi_\lambda(Z_{n})-\phi_\lambda(Z_{n+1})\;|\; \mathcal{F}_n]=\\
(e^{\lambda
Z_n}-1)\phi_\lambda(Z_{n})-\phi'_\lambda(Z_{n})\E[Z_{n+1}-Z_n\;|\;
\mathcal{F}_n]-\frac{1}{2}\phi''_\lambda(Z_{n})\E[(Z_{n+1}-Z_n)^2
\;|\; \mathcal{F}_n]-\E[\theta_n\;|\; \mathcal{F}_n].
\end{multline*}
In view of (a) and (b) of Lemma \ref{ART5_espzn+1-zn} and recalling the
differential equation (\ref{ART5_equadiff}) satisfied by $\phi_\lambda$,
the r.h.s. of the previous equality may be rewritten
\begin{equation*}
(e^{\lambda Z_n}-1-\lambda
Z_n)\phi_\lambda(Z_{n})-\phi'_\lambda(Z_{n})f_1(Z_n)-\phi''_\lambda(Z_{n})f_2(Z_n)-\E[\theta_n\;|\;
\mathcal{F}_n].
\end{equation*}
On the other hand, in view of (\ref{ART5_defWn}) and (\ref{ART5_defmu}), we have
\begin{equation}\label{ART5_defmusec}
\mu(n) = e^{-\lambda \sum_{k=0}^{n} Z_k}\E[e^{\lambda
Z_n}\phi_\lambda(Z_{n})-\phi_\lambda(Z_{n+1})\;|\; \mathcal{F}_n].
\end{equation}
Thus, for each $n\in[1,\sigma-1]$, we may decompose $\mu(n)$ in the form
\begin{equation}\label{ART5_decompnu}
\mu(n)=\mu_{1}(n)+\mu_{2}(n)+\mu_{3}(n)+\mu_{4}(n),
\end{equation}
where
\begin{eqnarray*}
\mu_{1}(n)&\defeg&e^{-\lambda
\sum_{k=0}^{n} Z_k}(e^{\lambda Z_n}-1-\lambda Z_n)\phi_\lambda(Z_n)\\
\mu_{2}(n)&\defeg&-e^{-\lambda
\sum_{k=0}^{n} Z_k}\phi'_\lambda(Z_{n})f_1(Z_n)\\
\mu_{3}(n)&\defeg&-e^{-\lambda
\sum_{k=0}^{n} Z_k}\phi''_\lambda(Z_{n})f_2(Z_n)\\
\mu_{4}(n)&\defeg&-e^{-\lambda
\sum_{k=0}^{n} Z_k}\E[\theta_n\;|\; \mathcal{F}_n].
\end{eqnarray*}
In particular, we can rewrite (\ref{ART5_alienZmu}) in the form (we
have to treat $\mu(0)$ separately since (\ref{ART5_deftheta}) does
not hold for $n=0$)
\begin{equation}\label{ART5_alienZmu2}
\E[1-e^{-\lambda \sum_{k=0}^{\sigma-1}
Z_k}]=\frac{1}{2^{\nu-1}\Gamma(\nu)}\left(\E\big[\mu(0)\big]
+\sum_{i=1}^4\E\Big[\sum_{n=1}^{\sigma-1} \mu_i(n)\Big]\right).
\end{equation}
We now state the main estimates:
\begin{lemm}\label{ART5_equivmun}There exist $\varepsilon>0$ and eight
finite constants $(C_i,C_i',i=0,2,3,4)$  such that, as $\lambda$
tends to $0^+$,
\begin{enumerate}
\item[\textup{(a)}]
$\E\left[\mu(0)\right]=\left\{\begin{array}{ll}
C_0\lambda^{\nu}+\mathcal{O}(\lambda)&\;\hspace*{1.4cm}\hbox{ if }\; \alpha\in(0,1)\\
C_0\lambda\log
\lambda+C_0'\lambda+o(\lambda)&\;\hspace*{1.4cm}\hbox{ if }\;
\alpha=1,\end{array}\right.$
\item[\textup{(b)}] $\E\left[\sum_{n=1}^{\sigma-1}
\mu_{1}(n)\right]= o(\lambda)\hspace*{4.3cm}\hbox{ for }
\;\alpha\in(0,1],$
\item[\textup{(c)}] $\E\left[\sum_{n=1}^{\sigma-1}
\mu_{2}(n)\right]=\left\{\begin{array}{ll}
C_2\lambda^{\nu}+o(\lambda^{\nu+\varepsilon})&\;\hbox{ if } \;\alpha\in(0,1)\\
C_2\lambda\log \lambda+C_2'\lambda+o(\lambda)&\;\hbox{ if }\;
\alpha=1,\end{array}\right.$
\item[\textup{(d)}] $\E\left[\sum_{n=1}^{\sigma-1}
\mu_{3}(n)\right]=\left\{\begin{array}{ll}
C_3\lambda^{\nu}+o(\lambda^{\nu+\varepsilon})&\;\hbox{ if } \;\alpha\in(0,1)\\
C_3\lambda\log \lambda+C_3'\lambda+o(\lambda)&\;\hbox{ if }\;
\alpha=1,\end{array}\right.$
\item[\textup{(e)}] $\E\left[\sum_{n=1}^{\sigma-1}
\mu_{4}(n)\right]=\left\{\begin{array}{ll}
C_4\lambda^{\nu}+o(\lambda^{\nu+\varepsilon})&\;\hspace*{1.4cm}\hbox{ if } \;\alpha\in(0,1)\\
C_4'\lambda+o(\lambda)&\;\hspace*{1.4cm}\hbox{ if }\;
\alpha=1.\end{array}\right.$
\end{enumerate}
\end{lemm}
Let us for the time being postpone the long and technical proof of
these estimates until section \ref{ART5_sectionpreuvemun} and complete
the proof of Proposition \ref{ART5_progeny}. In view of
(\ref{ART5_alienZmu2}), using the previous lemma, we deduce that there
exist some constants $C, C'$ such that
\begin{equation}\label{ART5_equivtotalxi}
\E\left[1-e^{-\lambda \sum_{k=0}^{\sigma-1}
Z_k}\right]=\left\{\begin{array}{ll}
C\lambda^{\nu}+o(\lambda^{\nu+\varepsilon})&\;\hbox{ if $\alpha\in(0,1)$,}\\
C\lambda \log \lambda+ C'\lambda+o(\lambda)&\;\hbox{ if $\alpha=1$.}\\
\end{array}\right.
\end{equation}
with
\begin{equation*}
C \defeg \left\{
\begin{array}{ll}
2^{1-\nu}\Gamma(\nu)^{-1}(C_0 + C_2 + C_3 + C_4)&\hbox{ when
$\alpha<1$,}\\
2^{1-\nu}\Gamma(\nu)^{-1}(C_0 + C_2 + C_3)&\hbox{ when
$\alpha=1$.}\\
\end{array}\right.
\end{equation*}
It simply remains to check that the constant $C$ is not zero.
Indeed, suppose that $C=0$. We first assume $\alpha=1$. Then, from
(\ref{ART5_equivtotalxi}),
\begin{equation*}
\E\left[1-e^{-\lambda \sum_{k=0}^{\sigma-1} Z_k}\right] =
C'\lambda+o(\lambda)
\end{equation*}
which implies $\E[\sum_{k=0}^{\sigma-1}Z_k]<\infty$ and contradicts
Corollary \ref{ART5_CorIntegZinfty}. Similarly, when $\alpha\in (0,1)$
and $C=0$, we get from (\ref{ART5_equivtotalxi}),
\begin{equation*}
\E\left[1-e^{-\lambda \sum_{k=0}^{\sigma-1}
Z_k}\right]=o(\lambda^{\nu+\varepsilon}).
\end{equation*}
This implies, for any $0<\varepsilon'<\varepsilon$, that
\begin{equation*} \E\left[(\sum_{n=0}^{\sigma-1} Z_n)^{\nu +
\varepsilon'}\right]<\infty
\end{equation*}
which contradicts Lemma \ref{ART5_ordreprogeny}. Therefore, $C$ cannot be
zero and the proposition is proved.
\end{proof}

\section{Technical estimates}\label{ART5_sectiontechnical}
\subsection{Some properties of modified Bessel functions}\label{ART5_sectionchoixphi}
We now recall some properties of modified Bessel functions. All the
results cited here may be found in \cite{AbramowitzStegun64}
(section 9.6) or \cite{Lebedev72} (section 5.7). For $\eta\in\R$,
the modified Bessel function of the first kind $I_\eta$ is defined
by
\begin{equation*}
I_\eta(x) \defeg
\left(\frac{x}{2}\right)^\eta\sum_{k=0}^{\infty}\frac{(x/2)^{2k}}{\Gamma(k+1)\Gamma(k+1+\eta)}
\end{equation*}
and the modified Bessel function of the second kind $K_\eta$ is given by the formula
\begin{equation*}
K_\eta(x) \defeg \left\{
\begin{array}{ll}
\frac{\pi}{2}\frac{I_{-\eta}(x)-I_\eta(x)}{\sin\pi\eta}&\hbox{ for
$\eta\in \R-\Z$,}\\
\lim_{\eta'\to \eta}K_{\eta'}(x)&\hbox{ for $\eta \in\Z$.}
\end{array}\right.
\end{equation*}
We are particularly interested in
\begin{equation*}
F_\eta(x) \defeg x^\eta K_\eta(x)\quad\hbox{for $x>0$.}
\end{equation*}
Thus, the function $\phi_\lambda$ defined in (\ref{ART5_defphilanbda1}) may be expressed in the form
\begin{equation}\label{ART5_defphi2}
\phi_\lambda(x)=F_\nu(\sqrt{\lambda}x).
\end{equation}
\begin{fact} For $\eta\geq 0$, the function $F_\eta$ is analytic,
positive and strictly decreasing on $(0,\infty)$. Moreover
\begin{enumerate}
\item \underline{Behaviour  at $0$}:
\begin{itemize}
\item[(a)] If $\eta>0$, the function $F_\eta$ is defined by continuity at 0 with $F_\eta(0)= 2^{\eta-1}\Gamma(\eta)$.
\item[(b)] If $\eta=0$, then  $F_0(x) = -\log x+\log 2-\gamma+o(1)$ as $x\to 0^+$ where $\gamma$ denotes Euler's constant.
\end{itemize}
\item \underline{Behaviour at infinity}:
$$F_\eta(x) \underset{x\to\infty}{\sim}\sqrt{\frac{\pi}{2x}}e^{-x}.$$
In particular, for every $\eta>0$, there exists $c_\eta\in \R$ such
that, for all $x\geq 0$,
\begin{equation}\label{ART5_majorationF}F_\eta(x)\le c_\eta
e^{-x}.\end{equation}
\item \underline{Formula for the derivative}:
\begin{equation}\label{ART5_deriveeF}
F'_\eta(x) = - x^{2\eta-1}F_{1-\eta}(x).
\end{equation}
In particular, $F_\eta$ solves the differential equation
\begin{equation*}
x F''_\eta(x)-(2\eta-1)F'_\eta(x)-xF_\eta(x)=0.
\end{equation*}
\end{enumerate}
\end{fact}
\noindent Concerning the function $\phi_\lambda$, in view of
(\ref{ART5_defphi2}), we deduce
\begin{fact}\label{ART5_reffact2} For each $\lambda>0$, the function $\phi_\lambda$ is
 analytic, positive and strictly decreasing on
$(0,\infty)$. Moreover
\begin{enumerate}
\item[\textup{(a)}] $\phi_\lambda$ is continuous and differentiable at $0$ with $\phi_\lambda(0) = 2^{\nu-1}\Gamma(\nu)$ and $\phi'_\lambda(0)=0$.
\item[\textup{(b)}] For $x>0$, we have
\begin{eqnarray*}
\phi'_\lambda(x)&=&-\lambda^{\nu} x^{\alpha} F_{1-\nu}(\sqrt{\lambda}x),\\
\phi''_\lambda(x)&=&\lambda F_{\nu}(\sqrt{\lambda}x)-\alpha \lambda^{\nu} x^{\alpha-1} F_{1-\nu}(\sqrt{\lambda}x).
\end{eqnarray*}
In particular, $\phi_\lambda$ solves the differential equation
\begin{equation*}
-\lambda
x\phi_{\lambda}(x)-\alpha\phi'_\lambda(x)+x\phi''_{\lambda}(x)=0.
\end{equation*}
\end{enumerate}
\end{fact}

\subsection{Proof of Lemma \ref{ART5_equivmun}}\label{ART5_sectionpreuvemun}
The proof of Lemma \ref{ART5_equivmun} is long and tedious
but requires only elementary methods. We shall treat, in separate
subsections the assertions (a) - (e) when $\alpha<1$. We explain, in
a last subsection, how to deal with the case $\alpha=1$.

We will use the following result extensively throughout the proof of
Lemma \ref{ART5_equivmun}.
\begin{lemm}\label{ART5_lemmutil} There exists $\varepsilon >0$ such that
\begin{equation*}
\E\left[\sigma(1-e^{-\lambda\sum_{k=0}^{\sigma-1}Z_k})\right] =
o(\lambda^\varepsilon)\quad\hbox{ as $\lambda\to 0^+$.}
\end{equation*}
\end{lemm}
\begin{proof} Let $\beta < \alpha \leq 1$, the function $x\to x^{\beta}$
is concave, thus
\begin{equation*}
\E\left[(\sum_{k=0}^{\sigma-1} Z_k)^\beta\right] \leq
\E\left[\sum_{k=0}^{\sigma-1} Z_k^\beta\right] \defeg c_1 < \infty,
\end{equation*}
where we used Corollary \ref{ART5_CorIntegZinfty} to conclude on the
finiteness of $c_1$. From Markov's inequality, we deduce that
$\P\left\{\sum_{k=0}^{\sigma-1} Z_k > x \right\} \leq
\frac{c_1}{x^{\beta}}$ for all $x\geq 0$. Therefore,
\begin{equation*}
\E\left[1-e^{-\lambda\sum_{k=0}^{\sigma-1}Z_k}\right] \leq
(1-e^{-\lambda x}) + \P\bigg\{\sum_{k=0}^{\sigma-1}Z_k > x\bigg\}
\leq \lambda x + \frac{c_1}{x^\beta}.
\end{equation*}
Choosing $x=\lambda^{-\frac{1}{\beta+1}}$ and setting $\beta' \defeg
\frac{\beta}{\beta+1}$, we deduce
\begin{equation*}
\E\left[1-e^{-\lambda\sum_{k=0}^{\sigma-1}Z_k}\right] \leq
(1+c_1)\lambda^{\beta'}.
\end{equation*}
According to Corollary \ref{ART5_approxsigma}, for $\delta<\alpha$, we
have $\E[\sigma^{1+\delta}]<\infty$, so Hölder's inequality gives
\begin{eqnarray*}
\E\left[\sigma(1-e^{-\lambda \sum_{k=0}^{\sigma-1}
Z_k})\right]&\leq&
\E[\sigma^{1+\delta}]^{\frac{1}{1+\delta}}\E\left[(1-e^{-\lambda
\sum_{k=0}^{\sigma-1}
Z_k})^{\frac{1+\delta}{\delta}}\right]^{\frac{\delta}{1+\delta}}\\
&\leq&
\E[\sigma^{1+\delta}]^{\frac{1}{1+\delta}}\E\left[1-e^{-\lambda
\sum_{k=0}^{\sigma-1} Z_k}\right]^{\frac{\delta}{1+\delta}} \,\leq\,
c_2 \lambda^{\frac{\beta'\delta}{1+\delta}},
\end{eqnarray*}
which completes the proof of the lemma.
\end{proof}

\subsubsection{Proof of (a) of Lemma \ref{ART5_equivmun} when
$\alpha<1$}\label{ART5_preuvea} Using the expression of $\mu(0)$ given by
(\ref{ART5_defmusec}) and the relation (\ref{ART5_deriveeF}) between of
$F'_{\nu}$ and $F_{1-\nu}$, we have
$$\E[\mu(0)]=\E[F_\nu(0)-F_\nu(\sqrt{\lambda}Z_1)]=-\E\left[\int_0^{\sqrt{\lambda}Z_1}F_\nu'(x)dx\right]
=\lambda^\nu\E\left[\int_0^{Z_1}y^\alpha
F_{1-\nu}(\sqrt{\lambda}y)dy\right].$$ Thus, using the dominated
convergence theorem,
\begin{equation*}
\lim_{\lambda\rightarrow
0}\frac{1}{\lambda^\nu}E[\mu(0)]=\E\left[\int_0^{Z_1}y^\alpha
F_{1-\nu}(0)dy\right]=\frac{F_{1-\nu}(0)}{1+\alpha}\E[Z_1^{1+\alpha}]
\defeg C_0 < \infty.
\end{equation*}
Furthermore, using again (\ref{ART5_deriveeF}), we get
\begin{eqnarray*}
\Big|\frac{1}{\lambda^\nu}E[\mu(0)]- C_0\Big|&=&\E\left[\int_0^{Z_1}y^\alpha\left(F_{1-\nu}(0)-
F_{1-\nu}(\sqrt{\lambda}y)\right)dy\right]\\
&=&\E\left[\int_0^{Z_1}y^\alpha
\int_0^{\sqrt{\lambda}y}x^{-\alpha}F_{\nu}(x)dx
dy\right]\\
&\le &
\frac{||F_\nu||_\infty}{1-\alpha}\lambda^{\frac{1-\alpha}{2}}\E\left[\int_0^{Z_1}y
dy\right] \,=\,\frac{||F_\nu||_\infty
\E[Z_1^2]}{2(1-\alpha)}\lambda^{\frac{1-\alpha}{2}}.
\end{eqnarray*}
Therefore, we obtain
$$E[\mu(0)]=C_0\lambda^{\nu}+\mathcal{O}(\lambda)$$
which proves (a) of Lemma \ref{ART5_equivmun}.

\subsubsection{Proof of (b) of Lemma \ref{ART5_equivmun} when $\alpha<1$}\label{ART5_debuttechnique}
Recall that
\begin{equation*}
\mu_{1}(n)=e^{-\lambda \sum_{k=0}^{n} Z_k}(e^{\lambda Z_n}-1-\lambda
Z_n)\phi_\lambda(Z_n)=e^{-\lambda \sum_{k=0}^{n} Z_k}(e^{\lambda
Z_n}-1-\lambda Z_n)F_\nu(\sqrt{\lambda}Z_{n}).
\end{equation*}
Thus,  $\mu_{1}(n)$ is almost surely positive and
\begin{equation*}
\mu_{1}(n) \leq (1-e^{-\lambda Z_n}-\lambda Z_n e^{-\lambda
Z_n})F_\nu(\sqrt{\lambda}Z_{n}).
\end{equation*}
Moreover, for any $y>0$, we have $1-e^{-y}- y e^{-y} \leq
\min(1,y^2)$, thus
\begin{eqnarray*}
\mu_{1}(n) &\leq& (1-e^{-\lambda Z_n}-\lambda Z_n e^{-\lambda
Z_n})F_\nu(\sqrt{\lambda}Z_{n}) \left( \Indic_{\{Z_n>\frac{-2\log
\lambda}{\sqrt{\lambda}}\}}+\Indic_{\{Z_n\le\frac{-2\log\lambda}{\sqrt{\lambda}}\}}\right)\\
&\leq &F_{\nu}(\sqrt{\lambda} Z_n)\Indic_{\{Z_n>\frac{-2\log
\lambda}{\sqrt{\lambda}}\}} \,+\, ||F_{\nu}||_{\infty}\lambda^2
Z_n^2\Indic_{\{Z_n\le\frac{-2\log\lambda}{\sqrt{\lambda}}\}}\\
&\leq& F_{\nu}(-2\log \lambda) \,+\, ||F_{\nu}||_{\infty}\lambda^2
Z_n^2\Indic_{\{Z_n\le\frac{-2\log\lambda}{\sqrt{\lambda}}\}},
\end{eqnarray*}
where we used the fact that $F_\nu$ is decreasing for the last
inequality. In view of (\ref{ART5_majorationF}), we also have
$F_{\nu}(-2\log\lambda)\le c_\nu \lambda^2$ and therefore
\begin{equation}\label{ART5_xip11}
\E\left[\sum_{n=1}^{\sigma-1} \mu_{1}(n)\right] \leq \lambda^2
c_\nu\E[\sigma] + \lambda^2
||F_{\nu}||_{\infty}\E\left[\sum_{n=1}^{\sigma-1}
Z_n^2\Indic_{\{Z_n\le\frac{-2\log \lambda}{\sqrt{\lambda}}\}}\right].
\end{equation}
On the one hand, according to (\ref{ART5_eqexc}), we have
\begin{equation}\label{ART5_xipp1}
\E\left[\sum_{n=1}^{\sigma-1} Z_n^2\Indic_{\{Z_n\le\frac{-2\log
\lambda}{\sqrt{\lambda}}\}}\right]=\E\left[
Z_\infty^2\Indic_{\{Z_\infty\le\frac{-2\log
\lambda}{\sqrt{\lambda}}\}}\right]\E[\sigma].
\end{equation}
On the other hand, Proposition \ref{ART5_proptailZinfty} states that
$\P(Z_\infty\geq x)\sim \frac{C}{x^\alpha}$ as $x$ tends to
infinity, thus
\begin{equation*}
\E\left[ Z_\infty^2\Indic_{\{Z_\infty\le
x\}}\right]\underset{x\to\infty}{\sim}
2\sum_{k=1}^{x}k\P(Z_\infty\ge k)\underset{x\to\infty}{\sim}
\frac{2C}{2-\alpha}x^{2-\alpha}.
\end{equation*}
This estimate and (\ref{ART5_xipp1}) yield
\begin{equation}\label{ART5_xi12}\lambda^2\E\left[\sum_{n=1}^{\sigma-1}
Z_n^2\Indic_{\{Z_n\le\frac{-2\log
\lambda}{\sqrt{\lambda}}\}}\right]\underset{\lambda\to 0^+}{\sim}
\Cste{3} \lambda^{1+\frac{\alpha}{2}}|\log \lambda|^{2-\alpha}.
\end{equation}
Combining (\ref{ART5_xip11}) and (\ref{ART5_xi12}), we finally obtain
\begin{equation*}
\E\left[\sum_{n=1}^{\sigma-1} \mu_{1}(n)\right]=o(\lambda),
\end{equation*}
which proves (b) of Lemma \ref{ART5_equivmun}.

\subsubsection{Proof of (c) of Lemma \ref{ART5_equivmun} when $\alpha<1$}
Recall that
\begin{equation*}
\mu_{2}(n)=-e^{-\lambda \sum_{k=0}^{n}
Z_k}\phi'_\lambda(Z_{n})f_1(Z_n)=\lambda^{\nu} Z_n^{\alpha}
F_{1-\nu}(\sqrt{\lambda}Z_n)f_1(Z_n)e^{-\lambda \sum_{k=0}^{n} Z_k}.
\end{equation*}
Since $f_1(x)= 0$ for $x\geq M-1$ (\emph{c.f.} Lemma
\ref{ART5_espzn+1-zn}),  the quantity $|\mu_{2}(n)|/\lambda^\nu$ is
smaller than $M^\alpha ||f_1||_\infty ||F_{1-\nu}||_\infty$. Thus,
using the dominated convergence theorem, we get
\begin{equation*}
\lim_{\lambda\rightarrow
0}\frac{1}{\lambda^\nu}\E\left[\sum_{n=1}^{\sigma-1}
\mu_{2}(n)\right]=\E\left[\sum_{n=1}^{\sigma-1}Z_n^{\alpha}
F_{1-\nu}(0)f_1(Z_n)\right] \defeg C_2\in \R.
\end{equation*}
 It remains to prove that, for
$\varepsilon>0$ small enough, as $\lambda\to 0^+$
\begin{equation}\label{ART5_inter1}
\Big|\frac{1}{\lambda^\nu}\E\left[\sum_{n=1}^{\sigma-1}
\mu_{2}(n)\right] - C_2 \Big| = o(\lambda^{\varepsilon}).
\end{equation}
We can rewrite the l.h.s. of (\ref{ART5_inter1}) in the form
\begin{multline}\label{ART5_inter2}
\Big|\E\left[\sum_{n=1}^{\sigma-1}Z_n^{\alpha}f_1(Z_n)
(F_{1-\nu}(0)-F_{1-\nu}(\sqrt{\lambda}Z_n))\right] \\+
\E\left[\sum_{n=1}^{\sigma-1}Z_n^{\alpha}f_1(Z_n)
F_{1-\nu}(\sqrt{\lambda}Z_n)(1-e^{-\lambda \sum_{k=0}^{n}
Z_k})\right]\Big|.
\end{multline}
On the one hand, the first term is bounded by
\begin{eqnarray*}\E\left[\sum_{n=1}^{\sigma-1}Z_n^{\alpha}|f_1(Z_n)|
(F_{1-\nu}(0)-F_{1-\nu}(\sqrt{\lambda}Z_n))\right]&\le&
M^\alpha||f_1||_{\infty}\E[\sigma]\int_{0}^{\sqrt{\lambda}M}|F'_{1-\nu}(x)|dx\\
&\le & M^\alpha||f_1||_{\infty}\E[\sigma] ||F_\nu||_{\infty}
\int_{0}^{\sqrt{\lambda}M}x^{1-2\nu}dx \\
&\le& \Cste{4} \lambda^{1-\nu},
\end{eqnarray*}
where we used formula (\ref{ART5_deriveeF}) for the expression of
$F'_{1-\nu}$ for the second inequality. On the other hand the second
term of (\ref{ART5_inter2}) is bounded by
\begin{eqnarray*}
\E\!\left[\sum_{n=1}^{\sigma-1}Z_n^{\alpha}|f_1(Z_n)|
F_{1-\nu}(\sqrt{\lambda}Z_n)(1\!-\!e^{-\lambda \sum_{k=0}^{n}
Z_k})\right]\!\!\!\!&\leq&\!\!\!\! M^\alpha||f_1||_{\infty}
||F_{1-\nu}||_{\infty}\E[\sigma(1\!-\!e^{-\lambda
\sum_{k=0}^{\sigma-1}
Z_k})]\\
&\leq&\!\!\! \Cste{5}\lambda^{\varepsilon}
\end{eqnarray*}
where we used Lemma \ref{ART5_lemmutil} for the last inequality. Putting
the pieces together, we conclude that (\ref{ART5_inter1}) holds for
$\varepsilon>0$ small enough.

\subsubsection{Proof of (d) of Lemma \ref{ART5_equivmun} when $\alpha<1$} Recall
that  \begin{eqnarray*}\mu_{3}(n)&\!=\!&-e^{-\lambda \sum_{k=0}^{n}
Z_k}\phi''_\lambda(Z_{n})f_2(Z_n)\\&\!=\!&-e^{-\lambda
\sum_{k=0}^{n} Z_k}f_2(Z_n)\left(\lambda
F_{\nu}(\sqrt{\lambda}Z_n)+\alpha \lambda^{\nu} Z_n^{\alpha-1}
F_{1-\nu}(\sqrt{\lambda}Z_n)\right).
\end{eqnarray*}
Note that, since $\alpha\leq 1$, we have $Z_n^{\alpha-1}\leq 1$ when
$Z_n\neq 0$. The quantities $f_2(Z_n)$, $F_{\nu}(\sqrt{\lambda}Z_n)$
and $F_{1-\nu}(\sqrt{\lambda}Z_n))$ are also bounded, so we check,
using the dominated convergence theorem, that
\begin{equation*}
\lim_{\lambda\rightarrow
0}\frac{1}{\lambda^\nu}\E\left[\sum_{n=1}^{\sigma-1}
\mu_{3}(n)\right]=-\alpha\E\left[\sum_{n=1}^{\sigma-1}Z_n^{\alpha-1}
F_{1-\nu}(0)f_2(Z_n)\right] \defeg C_3\in \R.
\end{equation*}
Furthermore we have
\begin{eqnarray}
\label{ART5_tempo1}\frac{1}{\lambda^\nu}\E\left[\sum_{n=1}^{\sigma-1}
\mu_{3}(n)\right] - C_3&=&
-\lambda^{1-\nu}\E\left[\sum_{n=1}^{\sigma-1}e^{-\lambda
\sum_{k=0}^{n}
Z_k}f_2(Z_n) F_{\nu}(\sqrt{\lambda}Z_n)\right]\\
\nonumber&& \,+\,
\alpha\E\left[\sum_{n=1}^{\sigma-1}Z_n^{\alpha-1}f_2(Z_n)
\left(F_{1-\nu}(0)-F_{1-\nu}(\sqrt{\lambda}Z_n)\right)\right]\\
\nonumber&& \,+\,
\alpha\E\left[\sum_{n=1}^{\sigma-1}Z_n^{\alpha-1}f_2(Z_n)
F_{1-\nu}(\sqrt{\lambda}Z_n)\left(1-e^{-\lambda \sum_{k=0}^{n}
Z_k}\right)\right].
\end{eqnarray}
The first term is clearly bounded by $\Cste{6}\lambda^{1-\nu}$. We
turn our attention to the second term. In view of (\ref{ART5_deriveeF}),
we have
\begin{equation*}
F_{1-\nu}(0)-F_{1-\nu}(\sqrt{\lambda}Z_n) =
\int_{0}^{\sqrt{\lambda}Z_n}\! x^{1-2\nu}F_{\nu}(x)dx \leq
\frac{||F_\nu||_\infty}{2-2\nu} \lambda^{1-\nu}Z_n^{2-2\nu} =
\frac{||F_\nu||_\infty}{1-\alpha} \lambda^{1-\nu}Z_n^{1-\alpha},
\end{equation*}
where we used $2-2\nu = 1-\alpha$ for the last equality. Therefore,
\begin{eqnarray*}
\Big|\E\!\left[\sum_{n=1}^{\sigma-1}Z_n^{\alpha-1}f_2(Z_n)
(F_{1-\nu}(0)-F_{1-\nu}(\sqrt{\lambda}Z_n))\right]\Big|& \leq &
\frac{||F_\nu||_\infty||f_2||_\infty}{1-\alpha} \lambda^{1-\nu}
\E\!\left[\sum_{n=1}^{\sigma-1}1\right]\\
& \leq & \frac{||F_\nu||_\infty||f_2||_\infty\E[\sigma]}{1-\alpha}
\lambda^{1-\nu}.
\end{eqnarray*}
As for the third term of (\ref{ART5_tempo1}), with the help of Lemma
\ref{ART5_lemmutil}, we find
\begin{eqnarray*}
\Big|\E\!\left[\sum_{n=1}^{\sigma-1}Z_n^{\alpha-1}f_2(Z_n)
F_{1-\nu}(\sqrt{\lambda}Z_n)(1\!-\!e^{-\lambda \sum_{k=0}^{n}
Z_k})\right]\!\Big|\!\!\!&\leq&\!\!\!
||f_2||_{\infty}||\!F_{1-\nu}||_{\infty}\E\!\left[\sigma(1\!-\!e^{-\lambda
\sum_{k=0}^{\sigma-1} Z_k})\right]\\
\!\!&\leq&\!\!\Cste{7}\lambda^{\varepsilon}.
\end{eqnarray*}
Putting the pieces together, we conclude that
\begin{equation*}
\E\left[\sum_{n=1}^{\sigma-1}
\mu_{3}(n)\right]=C_3\lambda^{\nu}+o(\lambda^{\nu+\varepsilon}).
\end{equation*}

\subsubsection{Proof of (e) of Lemma \ref{ART5_equivmun} when $\alpha<1$}
Recall that
\begin{equation}\label{ART5_defmu42}
\mu_{4}(n)=-e^{-\lambda \sum_{k=0}^{n} Z_k}\E[\theta_n\;|\;
\mathcal{F}_n].
\end{equation}
This term is clearly the most difficult to deal with. We first need
the next lemma stating that  $Z_{n+1}$ cannot be too "far" from
$Z_n$.
\begin{lemm}\label{ART5_largedeviation}
There exist  two constants $K_1,K_2>0$ such that for all $n\ge 0$,
\begin{enumerate}
\item[\textup{(a)}] $\P( Z_{n+1}\le \frac{1}{2}Z_n\;|\;\mathcal{F}_n)\le K_1e^{-K_2Z_n}$,
\item[\textup{(b)}] $\P(Z_{n+1}\ge 2 Z_n\;|\;\mathcal{F}_n)\le K_1e^{-K_2Z_n}$.
\end{enumerate}
\end{lemm}
\begin{proof}
This lemma follows from large deviation estimates. Indeed, with the
notation of section \ref{ART5_section2}, in view of Cramer's theorem, we
have, for any $j\geq M-1$,
\begin{eqnarray*}\P\Big\{ Z_{n+1}\le
\frac{1}{2}Z_n\;|Z_n=j\Big\}&=&
\P\Big\{A_{M-1}+\xi_1+\ldots+\xi_{j-M+1}\le \frac{j}{2}\Big\}\\
&\leq&\P\Big\{\xi_1+\ldots+\xi_{j-M+1}\le \frac{j}{2}\Big\} \,\leq\,
K_1e^{-K_2j},
\end{eqnarray*}
where we used the fact that $(\xi_i)$ is a sequence of i.i.d
geometric random variables with mean $1$. Similarly, recalling that
$A_{M-1}$ admits exponential moments of order $\beta < 2$, we also
deduce, for $j\geq M-1$, with possibly extended values of $K_1$
and $K_2$, that
\begin{multline*}\P\Big\{ Z_{n+1}\ge
2Z_n\;|Z_n=j\Big\}=
\P\Big\{A_{M-1}+\xi_1+\ldots+\xi_{j-M+1}\ge 2j\Big\}\\
\leq \P\Big\{A_{M-1}\ge \frac{j}{2}\Big\} +
\P\Big\{\xi_1+\ldots+\xi_{j-M+1}\ge \frac{3j}{2}\Big\} \leq
K_1e^{-K_2j}.
\end{multline*}
\end{proof}
Throughout this section, we use the notation, for $t\in[0,1]$ and
$n\in\N$,
\begin{equation*}
V_{n,t} \defeg Z_n+t(Z_{n+1}-Z_n).
\end{equation*}
In particular $V_{n,t}\in [Z_n,Z_{n+1}]$ (with the convention that
for $a>b$,  $[a,b]$ means $[b,a]$). With this notation, we can
rewrite the expression of $\theta_n$ given in (\ref{ART5_deftheta})
in the form
\begin{equation*}
\theta_n=(Z_{n+1}-Z_n)^2\int_{0}^{1}(1-t)\Big(
\phi_{\lambda}''(V_{n,t})-\phi_{\lambda}''(Z_n)\Big)dt.
\end{equation*}
Therefore, using the expression of $\phi_\lambda'$ and
$\phi_\lambda''$ stated in Fact (\ref{ART5_reffact2}), we get
\begin{equation}\label{ART5_thetaI1I2}
\E[\theta_n\;|\;\mathcal{F}_n] =
\int_{0}^{1}(1-t)(I_n^1(t)+I_n^2(t))dt,
\end{equation}
with
\begin{eqnarray*}
I_n^1(t)&\defeg&\lambda\E\left[(Z_{n+1}-Z_n)^2\Big(F_{\nu}(\sqrt{\lambda}V_{n,t})-
F_{\nu}(\sqrt{\lambda}Z_n)\Big)\;\Big|\;\mathcal{F}_n\right],\\
I_n^2(t)&\defeg&
-\alpha\lambda^{\nu}\E\left[(Z_{n+1}-Z_n)^2\left(V_{n,t}^{\alpha-1}F_{1-\nu}(\sqrt{\lambda}V_{n,t})-
Z_n^{\alpha-1}F_{1-\nu}(\sqrt{\lambda}Z_n)\right)\;\Big|\;\mathcal{F}_n\right].
\end{eqnarray*}
Recall that we want to estimate
\begin{multline*}
\E\!\left[\sum_{n=1}^{\sigma-1}\mu_4(n)\right] =
\E\!\left[\sum_{n=1}^{\sigma-1}e^{-\lambda \sum_{k=0}^{n}
Z_k}\int_{0}^{1}(1-t) I_n^1(t)dt\right] \\
+ \E\!\left[\sum_{n=1}^{\sigma-1}e^{-\lambda \sum_{k=0}^{n}
Z_k}\int_{0}^{1}(1-t) I_n^2(t)dt\right].
\end{multline*}
We deal with each term separately.

\bigskip \noindent\underline{\textbf{Dealing with $I^1$}}: We prove that the
contribution of this term is negligible, \emph{i.e.}
\begin{equation}\label{ART5_I1}
\Big|\E\left[\sum_{n=1}^{\sigma-1}e^{-\lambda \sum_{k=0}^{n} Z_k}
\int_{0}^{1}(1-t)I_n^1(t)dt\right]\Big|\leq \Cste{8}
\lambda^{\nu+\varepsilon}.
\end{equation}
To this end, we first notice that
\begin{eqnarray}
\nonumber|I_n^1(t)|&\leq&
\lambda^{\frac{3}{2}}\E\left[|Z_{n+1}-Z_n|^3
\max_{x\in[Z_n,Z_{n+1}]}|F'_{\nu}(\sqrt{\lambda}x)|\;\Big|\;\mathcal{F}_n
\right]\\
\nonumber&=&\lambda^{\frac{3}{2}}\E\left[|Z_{n+1}-Z_n|^3
\max_{x\in[Z_n,Z_{n+1}]}(\sqrt{\lambda}x)^\alpha
F_{1-\nu}(\sqrt{\lambda}x)\;\Big|\;\mathcal{F}_n \right]\\
\label{ART5_split1}&\le&
c_{1-\nu}\lambda^{\frac{3}{2}}\E\left[|Z_{n+1}-Z_n|^3
\max_{x\in[Z_n,Z_{n+1}]}(\sqrt{\lambda}x)^\alpha
e^{-\sqrt{\lambda}x}\;\Big|\;\mathcal{F}_n \right],
\end{eqnarray}
where we used (\ref{ART5_majorationF}) to find  $c_{1-\nu}$ such that
$F_{1-\nu}(x)\le c_{1-\nu} e^{-x}$. We now split (\ref{ART5_split1})
according to whether
\begin{equation*}
\hbox{(a) }\frac{1}{2}Z_n\le Z_{n+1}\le 2 Z_n\quad\hbox{ or
}\quad\hbox{(b) } Z_{n+1}<\frac{1}{2}Z_n \hbox{ or } Z_{n+1}>2 Z_n.
\end{equation*}
One the one hand, Lemma \ref{ART5_espzn+1-zn} states that
\begin{equation*}
\E\left[|Z_{n+1}-Z_n|^p\;|\;\mathcal{F}_n \right]\le D_p
Z_n^{\frac{p}{2}} \quad \mbox{ for all } p\in \N \mbox{ and }
Z_n\neq 0.
\end{equation*}
Hence, for $1\leq n\leq \sigma-1$, we get
\begin{multline}\label{ART5_split2}
\E\left[|Z_{n+1}-Z_n|^3
\max_{x\in[Z_n,Z_{n+1}]}(\sqrt{\lambda}x)^\alpha
e^{-\sqrt{\lambda}x}\Indic_{\{\frac{1}{2}Z_n\le Z_{n+1}\le 2
Z_n\}}\;\Big|\;\mathcal{F}_n \right]\\
\begin{aligned}
&\le  \E\left[|Z_{n+1}-Z_n|^3
\max_{x\in[\frac{1}{2}Z_n,2Z_{n}]}(\sqrt{\lambda}x)^\alpha
e^{-\sqrt{\lambda}x}\Indic_{\{\frac{1}{2}Z_n\le Z_{n+1}\le 2
Z_n\}}\;\Big|\;\mathcal{F}_n \right]\\
& \le  \E\left[|Z_{n+1}-Z_n|^3 (2\sqrt{\lambda}Z_n)^\alpha
e^{-\frac{1}{2}\sqrt{\lambda}Z_n}\;\Big|\;\mathcal{F}_n \right]\\
&\le \Cste{9} Z_n^{\frac{3}{2}}(\sqrt{\lambda}Z_n)^\alpha
e^{-\frac{1}{2}\sqrt{\lambda}Z_n}\\
&\le \Cste{9}\lambda^{\frac{3\alpha-6}{8}}
Z_n^{\frac{3\alpha}{4}}(\sqrt{\lambda}Z_n)^{\frac{6+\alpha}{4}}
e^{-\frac{1}{2}\sqrt{\lambda}Z_n}\\
&\le \Cste{10} \lambda^{\frac{3\alpha-6}{8}}
Z_n^{\frac{3\alpha}{4}},
\end{aligned}
\end{multline}
where we used the fact that the function $x^\frac{6+\alpha}{4}
e^{-\frac{x}{2}}$ is bounded on $\R_+$ for the last inequality. On
the other hand,
\begin{multline}\label{ART5_split3}
\E\left[|Z_{n+1}-Z_n|^3
\max_{x\in[Z_n,Z_{n+1}]}(\sqrt{\lambda}x)^\alpha
e^{-\sqrt{\lambda}x}\Indic_{\{Z_{n+1}<\frac{1}{2}Z_n \mbox{ or }
Z_{n+1}> 2
Z_n\}}\;\Big|\;\mathcal{F}_n \right]\\
\begin{aligned}
&\le  \E\left[|Z_{n+1}-Z_n|^3 \max_{x\ge 0}(\sqrt{\lambda}x)^\alpha
e^{-\sqrt{\lambda}x}\Indic_{\{Z_{n+1}<\frac{1}{2}Z_n \mbox{ or }
Z_{n+1}> 2
Z_n\}}\;\Big|\;\mathcal{F}_n \right]\\
&\le \Cste{11}\E\left[|Z_{n+1}-Z_n|^6\;|\;\mathcal{F}_n
\right]^{1/2} \P\left\{Z_{n+1}<\frac{1}{2}Z_n \mbox{ or } Z_{n+1}> 2
Z_n\;\Big|\;\mathcal{F}_n \right\}^{\frac{1}{2}}\\
&\le \Cste{12}Z_n^{\frac{3}{2}}e^{-\frac{K_2}{2}Z_n}\\
&\le \Cste{13}.
\end{aligned}
\end{multline}
Combining (\ref{ART5_split1}), (\ref{ART5_split2}) and (\ref{ART5_split3}), we get
\begin{equation*}
|I_n^1(t)|\;\leq\;
c_{1-\nu}\Cste{13}\lambda^{\frac{3}{2}}+c_{1-\nu}\Cste{10}
\lambda^{\frac{3\alpha+6}{8}} Z_n^{\frac{3\alpha}{4}} \;\leq\;
\Cste{14}\lambda^{\nu+\frac{2-\alpha}{8}}Z_n^{\frac{3\alpha}{4}}.
\end{equation*}
And therefore
\begin{equation*}
\Big|\E\left[\sum_{n=1}^{\sigma-1}e^{-\lambda \sum_{k=0}^{n} Z_k}
\int_{0}^{1}(1-t)I_n^1(t)dt\right]\Big|\leq
\Cste{14}\lambda^{\nu+\frac{2-\alpha}{8}}\E\left[\sum_{n=1}^{\sigma-1}Z_n^{\frac{3\alpha}{4}}\right].
\end{equation*}
Corollary \ref{ART5_CorIntegZinfty} states that
$\E[\sum_{n=1}^{\sigma-1}Z_n^{\frac{3\alpha}{4}}]$ is finite so the
proof of (\ref{ART5_I1}) is complete.

\bigskip
\noindent\underline{\textbf{Dealing with $I^2$}}: It remains to
prove that
\begin{equation}\label{ART5_toprove}
\E\left[\sum_{n=1}^{\sigma-1}e^{-\lambda \sum_{k=0}^{n} Z_k}
\int_{0}^{1}(1-t)I_n^2(t)dt\right] = C_4 \lambda^\nu +
o(\lambda^{\nu+\varepsilon}).
\end{equation}
To this end, we write
\begin{equation}\label{ART5_decompI2}
I_n^2(t)=-\alpha \lambda^{\nu}(J_n^1(t)+J_n^2(t)+J_n^3(t)),
\end{equation}
with
\begin{eqnarray*}
J_n^1(t)&\defeg&
\E\left[(Z_{n+1}-Z_n)^2(F_{1-\nu}(\sqrt{\lambda}V_{n,t}))-F_{1-\nu}(\sqrt{\lambda}Z_n))
 Z_n^{\alpha-1}\;|\;\mathcal{F}_n \right], \\
J_n^2(t)&\defeg&
\E\left[(Z_{n+1}-Z_n)^2(V_{n,t}^{\alpha-1}-Z_n^{\alpha-1})
(F_{1-\nu}(\sqrt{\lambda}V_{n,t})-
 F_{1-\nu}(0))\;|\;\mathcal{F}_n \right],\\
J_n^3(t)&\defeg&F_{1-\nu}(0)
\E\left[(Z_{n+1}-Z_n)^2(V_{n,t}^{\alpha-1}-
 Z_n^{\alpha-1})\;|\;\mathcal{F}_n \right].
\end{eqnarray*}
Again, we shall study each term separately. In view of
(\ref{ART5_toprove}) and (\ref{ART5_decompI2}), the proof of (e) of Lemma
\ref{ART5_equivmun}, when $\alpha<1$, will finally be complete once we
established the following three estimates:
\begin{eqnarray}
\label{ART5_J1proof}\E\left[\sum_{n=1}^{\sigma-1}e^{-\lambda
\sum_{k=0}^{n} Z_k}
\int_{0}^{1}(1-t)J_n^1(t)dt\right]&=&\mathcal{O}(
\lambda^{\frac{1-\alpha}{4}}),\\
\label{ART5_J2proof}\E\left[\sum_{n=1}^{\sigma-1}e^{-\lambda
\sum_{k=0}^{n} Z_k}
\int_{0}^{1}(1-t)J_n^2(t)dt\right]&=& o(\lambda^{\varepsilon}),\\
\label{ART5_J3proof}\E\left[\sum_{n=1}^{\sigma-1}e^{-\lambda
\sum_{k=0}^{n} Z_k} \int_{0}^{1}(1-t)J_n^3(t)dt\right]&=&
C+o(\lambda^{\varepsilon}).
\end{eqnarray}

\bigskip
\noindent\underline{\textbf{Proof of (\ref{ART5_J1proof})}}: Using a
technique similar to that used for $I^1$, we split $J^1$ into three
different terms according to whether
\begin{equation*}
\hbox{(a) }\frac{1}{2}Z_n\leq Z_{n+1}\qquad\hbox{(b) }1\le
Z_{n+1}<\frac{1}{2}Z_n\qquad\hbox{(c) }Z_{n+1}=0.
\end{equation*}
For the first case (a), we write, for $1\leq n\leq \sigma-1$,
recalling that $V_{n,t}\in [Z_n,Z_{n+1}]$,
\begin{multline}\label{ART5_inter4}
\Big|\E\left[(Z_{n+1}-Z_n)^2\left(F_{1-\nu}(\sqrt{\lambda}V_{n,t})-F_{1-\nu}(\sqrt{\lambda}Z_n)\right)
 Z_n^{\alpha-1}\Indic_{\{\frac{1}{2}Z_n\le Z_{n+1}\}}\;\Big|\;\mathcal{F}_n \right]\Big|\\
\begin{aligned}
&\le  \lambda^{\frac{1}{2}}\E\Big[|Z_{n+1}-Z_n|^3
Z_n^{\alpha-1}\max_{x\ge
\frac{1}{2}Z_n}|F'_{1-\nu}(\sqrt{\lambda}x)|\;\Big|\;\mathcal{F}_n \Big]\\
&= \lambda^{\frac{1}{2}}\E\left[|Z_{n+1}-Z_n|^3\;|\;\mathcal{F}_n
\right] Z_n^{\alpha-1}\max_{x\ge
\frac{1}{2}Z_n}\left((\sqrt{\lambda}x)^{-\alpha}F_\nu(\sqrt{\lambda}x)\right)\\
 &\le  \Cste{15}\lambda^{\frac{1}{2}}\E\left[|Z_{n+1}-Z_n|^3\;|\;\mathcal{F}_n
\right] Z_n^{\alpha-1} \max_{x\ge
\frac{1}{2}Z_n}\left((\sqrt{\lambda}x)^{-\alpha}
e^{-\sqrt{\lambda}x}\right)\\
& = \Cste{15}
\lambda^{\frac{1}{2}}\E\left[|Z_{n+1}-Z_n|^3\;|\;\mathcal{F}_n
\right] Z_n^{-1}(\frac{1}{2}\sqrt{\lambda})^{-\alpha}
e^{-\frac{1}{2}\sqrt{\lambda}Z_n}\\
&\le \Cste{16} Z_n^{\frac{1}{2}}\lambda^{\frac{1-\alpha}{2}}
e^{-\frac{1}{2}\sqrt{\lambda}Z_n}\\
& = \Cste{16}\lambda^{\frac{1-\alpha}{4}}
Z_n^{\frac{\alpha}{2}}\left((\sqrt{\lambda}Z_n)^{\frac{1-\alpha}{2}}
e^{-\frac{1}{2}\sqrt{\lambda}Z_n}\right)\\
&\le \Cste{17}\lambda^{\frac{1-\alpha}{4}} Z_n^{\frac{\alpha}{2}},
\end{aligned}
\end{multline}
where we used Lemma \ref{ART5_espzn+1-zn} to get an upper bound for
the conditional expectation.

\noindent For the second case (b), keeping in mind Lemma
\ref{ART5_largedeviation}, we get
\begin{multline}\label{ART5_inter5}
\E\Big[(Z_{n+1}-Z_n)^2\left(F_{1-\nu}(\sqrt{\lambda}V_{n,t})-F_{1-\nu}(\sqrt{\lambda}Z_n)\right)
 Z_n^{\alpha-1}\Indic_{\{1\le Z_{n+1}<\frac{1}{2}Z_n\}}\;\Big|\;\mathcal{F}_n \Big]\\
\begin{aligned}
&\le  \Cste{18}\lambda^{\frac{1}{2}}\E\left[|Z_{n+1}-Z_n|^3
Z_n^{\alpha-1} \Indic_{\{1\le
Z_{n+1}<\frac{1}{2}Z_n\}}\;|\;\mathcal{F}_n \right]\max_{x\ge
1}\left((\sqrt{\lambda}x)^{-\alpha}
e^{-\sqrt{\lambda}x}\right)\\
& \le \Cste{19} \lambda^{\frac{1}{2}}\E\left[ Z_n^{\alpha+2}
\Indic_{\{1\le Z_{n+1}<\frac{1}{2}Z_n\}}\;|\;\mathcal{F}_n \right]\lambda^{-\frac{\alpha}{2}}\\
&\le \Cste{19} \lambda^{\frac{1-\alpha}{2}} Z_n^{\alpha+2} \P\{
Z_{n+1}<\frac{1}{2}Z_n\;|\;\mathcal{F}_n \}\\
&\le \Cste{19}K_1 \lambda^{\frac{1-\alpha}{2}} Z_n^{\alpha+2}e^{-K_2Z_n}\\
 &\le \Cste{20}\lambda^{\frac{1-\alpha}{2}}.
\end{aligned}
\end{multline}
For the last case (c), we note that when $Z_{n+1}=0$, then $V_{n,t}
= (1-t)Z_n$, therefore
\begin{multline}\label{ART5_inter6}
\E\Big[(Z_{n+1}-Z_n)^2\left(F_{1-\nu}(\sqrt{\lambda}V_{n,t})-F_{1-\nu}(\sqrt{\lambda}Z_n)\right)
 Z_n^{\alpha-1}\Indic_{\{Z_{n+1}=0\}}\;\Big|\;\mathcal{F}_n \Big]\\
\begin{aligned}
&=Z_n^2(F_{1-\nu}(\sqrt{\lambda}(Z_n(1-t)))-F_{1-\nu}(\sqrt{\lambda}Z_n))
 Z_n^{\alpha-1}\P\{Z_{n+1}=0\;|\;\mathcal{F}_n \}\\
&\le \Cste{21}\lambda^{\frac{1}{2}}Z_n^{2+\alpha} e^{-K_2
Z_n}\max_{x\in[Z_n(1-t),Z_n]}
(\sqrt{\lambda}x)^{-\alpha} \\
 &\le \Cste{21}
\lambda^{\frac{1-\alpha}{2}}(1-t)^{-\alpha}Z_n^{2}
 e^{-K_2
Z_n}\\
&\leq \Cste{22} \lambda^{\frac{1-\alpha}{2}}(1-t)^{-\alpha}.
\end{aligned}
\end{multline}
Combining (\ref{ART5_inter4}), (\ref{ART5_inter5}) and (\ref{ART5_inter6}), we
deduce that, for $1\leq n\leq \sigma-1$,
\begin{equation*}
\int_0^1(1-t)|J_n^1(t)|dt\le \Cste{23}\lambda^{\frac{1-\alpha}{4}}
Z_n^{\frac{\alpha}{2}}.
\end{equation*}
Moreover, according to Corollary \ref{ART5_CorIntegZinfty}, we have
$\E\left[\sum_{n=1}^{\sigma-1}Z_n^{\frac{\alpha}{2}}\right]<\infty$,
therefore
\begin{equation}\label{ART5_I3}\Big|\E\left[\sum_{n=1}^{\sigma-1}e^{-\lambda
\sum_{k=0}^{n} Z_k} \int_{0}^{1}(1-t)J_n^1(t)dt\right]\Big|\le
\E\left[\sum_{n=1}^{\sigma-1}
\int_{0}^{1}(1-t)|J_n^1(t)|dt\right]\le \Cste{24}
\lambda^{\frac{1-\alpha}{4}}
\end{equation}
which yields (\ref{ART5_J1proof}).

\bigskip
\noindent\underline{\textbf{Proof of (\ref{ART5_J2proof})}}: We write
\begin{equation*}
 J_n^2(t) = \E[R_n(t)\;|\;\mathcal{F}_n]
\end{equation*}
with
\begin{equation*}
R_n(t) \defeg
(Z_{n+1}-Z_n)^2\left(V_{n,t}^{\alpha-1}-Z_n^{\alpha-1}\right)
\left(F_{1-\nu}(\sqrt{\lambda}V_{n,t})-
 F_{1-\nu}(0)\right).
\end{equation*}
Again, we split the expression of $J^2$ according to four cases:
\begin{eqnarray}
\nonumber J_n^2(t) &=&
\E[R_n(t)\Indic_{\{Z_{n+1}=0\}}\;|\;\mathcal{F}_n] +
\E[R_n(t)\Indic_{\{1\le
Z_{n+1}<\frac{1}{2}Z_n\}}\;|\;\mathcal{F}_n]\\
\label{ART5_J2decoupe}&&+ \E[R_n(t)\Indic_{\{\frac{1}{2}Z_n\le Z_{n+1}\le 2
Z_n\}}\;|\;\mathcal{F}_n] + \E[R_n(t)\Indic_{\{Z_{n+1}> 2
Z_n\}}\;|\;\mathcal{F}_n].
\end{eqnarray}
We do not detail the cases $Z_{n+1}=0$ and $1\le
Z_{n+1}<\frac{1}{2}Z_n$ which may be treated by the same method
used in (\ref{ART5_inter5}) and (\ref{ART5_inter6}) and yields
similar bounds which do not depend on $Z_n$:
\begin{eqnarray*}
\E[R_n(t)\Indic_{\{Z_{n+1}=0\}}\;|\;\mathcal{F}_n] &\leq& \Cste{25}\lambda^{\frac{1-\alpha}{2}}(1-t)^{-\alpha}\\
\E[R_n(t)\Indic_{\{1\le Z_{n+1}<\frac{1}{2}Z_n\}}\;|\;\mathcal{F}_n]
&\leq& \Cste{26} \lambda^{\frac{1-\alpha}{2}}.
\end{eqnarray*}
In particular, the combination of these two estimates gives:
\begin{equation}\label{ART5_combi1}\left|\E\left[\sum_{n=1}^{\sigma-1}e^{-\lambda
\sum_{k=0}^{n} Z_k} \int_{0}^{1}(1-t)\E[R_n(t)\Indic_{\{Z_{n+1}<
\frac{Z_n}{2}\}}\;|\;\mathcal{F}_n]dt\right]\right|\le \Cste{27}
\lambda^{\frac{1-\alpha}{2}}.
\end{equation}

\bigskip
\noindent In order to deal with the third term on the r.h.s. of
(\ref{ART5_J2decoupe}), we write
\begin{multline*}
|\E[R_n(t)\Indic_{\{\frac{1}{2}Z_n\le Z_{n+1}\le 2
Z_n\}}\;|\;\mathcal{F}_n]|\\
\begin{aligned}
&=\Big|\E\Big[(Z_{n+1}-Z_n)^2(V_{n,t}^{\alpha-1}-Z_n^{\alpha-1})
(F_{1-\nu}(\sqrt{\lambda}V_{n,t})-
F_{1-\nu}(0))\Indic_{\{\frac{1}{2}Z_n\le Z_{n+1}\le 2
Z_n\}}\;|\;\mathcal{F}_n
\Big]\Big|\\
&\le \Cste{28} \E\left[|Z_{n+1}-Z_n|^3 \max_{x\ge \frac{Z_n}{2}}
x^{\alpha-2}\int_0^{2\sqrt{\lambda}Z_n}|F'_{1-\nu}(y)|dy\;\Big|\;\mathcal{F}_n\right]\\
&\le \Cste{29} \E\left[|Z_{n+1}-Z_n|^3\;|\;\mathcal{F}_n\right]
\max_{x\ge \frac{Z_n}{2}}
x^{\alpha-2}\int_0^{2\sqrt{\lambda}Z_n}y^{-\alpha}dy\\
&\le \Cste{30}\lambda^{\frac{1-\alpha}{2}} Z_n^{\frac{1}{2}}.
\end{aligned}
\end{multline*}
According to Corollary \ref{ART5_CorIntegZinfty}, when
$\frac{1}{2}<\alpha<1$, we have
$\E\left[\sum_{n=1}^{\sigma-1}Z_n^{1/2}\right]<\infty$. In this
case, we get
\begin{equation}\label{ART5_upbo}\left|\E\left[\sum_{n=1}^{\sigma-1}e^{-\lambda
\sum_{k=0}^{n} Z_k}
\int_{0}^{1}(1-t)\E[R_n(t)\Indic_{\{\frac{1}{2}Z_n\le Z_{n+1}\le 2
Z_n\}}\;|\;\mathcal{F}_n]dt\right]\right|\le \Cste{31}
\lambda^{\frac{1-\alpha}{2}}.
\end{equation}
When $0 < \alpha \leq\frac{1}{2}$, the function
$x^\frac{2-3\alpha}{4}e^{-x}$ is bounded on $\R_+$, so
\begin{eqnarray*}
e^{-\lambda Z_n}\int_0^1(1-t)|\E[R_n(t)\Indic_{\{\frac{1}{2}Z_n\le
Z_{n+1}\le 2 Z_n\}}\;|\;\mathcal{F}_n]| dt&\le&
\Cste{30}\lambda^{\frac{\alpha}{4}}Z_n^{\frac{3\alpha}{4}} (\lambda
Z_n)^{\frac{2-3\alpha}{4}}e^{-\lambda Z_n}\\
&\le& \Cste{32}\lambda^{\frac{\alpha}{4}}Z_n^{\frac{3\alpha}{4}}.
\end{eqnarray*}
Therefore, when $\alpha\leq \frac{1}{2}$, the estimate
(\ref{ART5_upbo}) still holds by changing
$\lambda^{\frac{1-\alpha}{2}}$ to $\lambda^{\frac{\alpha}{4}}$.
Hence, for every $\alpha\in (0,1)$, we can find $\varepsilon>0$ such
that
\begin{equation}\label{ART5_combi2}\left|\E\left[\sum_{n=1}^{\sigma-1}e^{-\lambda
\sum_{k=0}^{n} Z_k}
\int_{0}^{1}(1-t)\E[R_n(t)\Indic_{\{\frac{1}{2}Z_n\le Z_{n+1}\le 2
Z_n\}}\;|\;\mathcal{F}_n]dt\right]\right|\le \Cste{33}
\lambda^{\varepsilon}.
\end{equation}

\bigskip
\noindent We now give the upper bound for the last term on the
r.h.s. of (\ref{ART5_J2decoupe}). We have
\begin{eqnarray*}
\E\Big[R_n(t)\Indic_{\{Z_{n+1}\ge 2 Z_n\}}\;\Big|\;\mathcal{F}_n\Big]
&=& \E\Big[R_n(t)\Indic_{\{2Z_n \leq Z_{n+1}\leq
\lambda^{-\frac{1}{4}}\}}\;\Big|\;\mathcal{F}_n\Big]\\
&& + \E\Big[R_n(t)\Indic_{\{Z_{n+1}>
\max(\lambda^{-\frac{1}{4}},2Z_n)\}}\;\Big|\;\mathcal{F}_n\Big].
\end{eqnarray*}
On the one hand, when $Z_n \neq 0$ and $Z_{n+1}\neq 0$, we have
$|V_{n,t}^{\alpha-1}-Z_n^{\alpha-1}| \leq 2$ thus, for $1\leq n\leq
\sigma-1$,
\begin{multline*}
\Big|\E\Big[R_n(t)\Indic_{\{2Z_n \leq Z_{n+1}\leq
\lambda^{-\frac{1}{4}}\}}\;\Big|\;\mathcal{F}_n\Big]\Big|\\
\begin{aligned}
&=
\Big|\E\left[(Z_{n+1}-Z_n)^2\Big(V_{n,t}^{\alpha-1}-Z_n^{\alpha-1}\Big)
\Big(F_{1-\nu}(\sqrt{\lambda}V_{n,t})-
 F_{1-\nu}(0)\Big)\Indic_{\{2Z_n<Z_{n+1}\le \lambda^{-\frac{1}{4}}\}}\;\Big|\;\mathcal{F}_n
 \right]\Big|\\
&\leq
2\E\Big[(Z_{n+1}-Z_n)^2\int_{0}^{\sqrt{\lambda}Z_{n+1}}x^{-\alpha}F_\nu(x)dx
\Indic_{\{2Z_n<Z_{n+1}\le \lambda^{-\frac{1}{4}}\}}\;|\;\mathcal{F}_n
 \Big]\\
&\leq
\Cste{34}\E\Big[(Z_{n+1}-Z_n)^2\int_{0}^{\lambda^{\frac{1}{4}}}x^{-\alpha}dx
\Indic_{\{Z_{n+1}> 2Z_n\}}\;|\;\mathcal{F}_n
\Big]\\
&\leq \Cste{35}\lambda^{\frac{1-\alpha}{4}}\E\Big[(Z_{n+1}-Z_n)^2
\Indic_{\{Z_{n+1}> 2Z_n\}}\;|\;\mathcal{F}_n
\Big]\\
&\leq
\Cste{35}\lambda^{\frac{1-\alpha}{4}}\E\Big[(Z_{n+1}-Z_n)^4\;|\;\mathcal{F}_n
\Big]^{\frac{1}{2}}\P\Big\{Z_{n+1}>
2Z_n\;|\;\mathcal{F}_n\Big\}^{\frac{1}{2}}\\
&\leq \Cste{36}\lambda^{\frac{1-\alpha}{4}},
\end{aligned}
\end{multline*}
where we used Lemma \ref{ART5_espzn+1-zn} and Lemma \ref{ART5_largedeviation}
for the last inequality. On the other hand,
\begin{multline*}
\E\Big[R_n(t)\Indic_{\{Z_{n+1}>
\max(\lambda^{-\frac{1}{4}},2Z_n)\}}\;\Big|\;\mathcal{F}_n\Big]\\
\begin{aligned}
&= \Big|\E\Big[(Z_{n+1}-Z_n)^2(V_{n,t}^{\alpha-1}-Z_n^{\alpha-1})
(F_{1-\nu}(\sqrt{\lambda}V_{n,t})-
 F_{1-\nu}(0))\Indic_{\{Z_{n+1}> \max(\lambda^{-\frac{1}{4}},2Z_n)\}}\;|\;\mathcal{F}_n
 \Big]\Big|\\
&\le 2||F_{1-\nu}||_\infty \E\Big[(Z_{n+1}-Z_n)^2 \Indic_{\{Z_{n+1}>
\max(\lambda^{-\frac{1}{4}},2Z_n)\}}\;|\;\mathcal{F}_n
 \Big]\\
&\le \Cste{37} \E\Big[(Z_{n+1}-Z_n)^4
\Indic_{\{Z_{n+1}>2Z_n\}}\;|\;\mathcal{F}_n
\Big]^{\frac{1}{2}}\P\{Z_{n+1}>\lambda^{-\frac{1}{4}}\;|\;\mathcal{F}_n
 \}^{\frac{1}{2}}\\
&\leq
\Cste{38}Z_ne^{-\frac{K_2}{4}Z_n}\P\{Z_{n+1}>\lambda^{-\frac{1}{4}}\;|\;\mathcal{F}_n
 \}^{\frac{1}{2}}\\
&\leq \Cste{38}Z_ne^{-\frac{K_2}{4}Z_n}\E[Z_{n+1}\;|\;\mathcal{F}_n
]^{\frac{1}{2}}\lambda^{\frac{1}{8}}\\
&\leq\Cste{39}\lambda^{\frac{1}{8}}.
\end{aligned}
\end{multline*}
These two bounds yield
\begin{equation}\label{ART5_combi3}\left|\E\left[\sum_{n=1}^{\sigma-1}e^{-\lambda
\sum_{k=0}^{n} Z_k}
\int_{0}^{1}(1-t)\E[R_n(t)\Indic_{\{\frac{1}{2}Z_n\le Z_{n+1} > 2
Z_n\}}\;|\;\mathcal{F}_n]dt\right]\right|\le \Cste{40}
\lambda^{\beta}
\end{equation}
with $\beta=\min(\frac{1-\alpha}{4},\frac{1}{8})$. Combining
(\ref{ART5_combi1}), (\ref{ART5_combi2}) and (\ref{ART5_combi3}), we finally obtain
(\ref{ART5_J2proof}).

\bigskip
\noindent\underline{\textbf{Proof of (\ref{ART5_J3proof})}}: Recall that

\begin{equation*}
J_n^3(t) \defeg F_{1-\nu}(0)
\E\left[(Z_{n+1}-Z_n)^2(V_{n,t}^{\alpha-1}-
 Z_n^{\alpha-1})\;|\;\mathcal{F}_n \right].
\end{equation*}
In particular, $J_n^3(t)$ does not depend on $\lambda$. We want to
show that there exist $C\in \R$ and $\varepsilon>0$ such that
\begin{equation}\label{ART5_I4} \E\left[\sum_{n=1}^{\sigma-1}e^{-\lambda \sum_{k=0}^{n}
Z_k} \int_{0}^{1}(1-t)J_n^2(t)dt\right]= C+
o(\lambda^{\varepsilon}).
\end{equation}
We must first check that
\begin{equation*} \E\left[\sum_{n=1}^{\sigma-1}
\int_{0}^{1}(1-t)|J_n^2(t)|dt\right]<\infty.
\end{equation*}
This may be done, using the same method as before by distinguishing
three cases:
\begin{equation*}
\hbox{(a) }Z_{n+1}\ge \frac{1}{2}Z_n\qquad\hbox{(b) }1\le
Z_{n+1}<\frac{1}{2}Z_n\qquad\hbox{(c) }Z_{n+1}=0.
\end{equation*}
Since the arguments are very similar to those provided above, we
feel free to skip the details. We find, for $1\leq n\leq \sigma-1$,
\begin{equation*}
\int_{0}^{1}(1-t)|J_n^2(t)|dt\;\le\;\Cste{41}Z_n^{\alpha-\frac{1}{2}}+\Cste{42}
\;\le\; \Cste{43}Z_n^{\frac{\alpha}{2}}.
\end{equation*}
Since
$\E\left[\sum_{n=1}^{\sigma-1}Z_n^{\frac{\alpha}{2}}\right]<\infty$,
with the help of the dominated convergence theorem, we get
\begin{equation*}\lim_{\lambda\rightarrow 0} \E\left[\sum_{n=1}^{\sigma-1}e^{-\lambda \sum_{k=0}^{n}
Z_k} \int_{0}^{1}(1-t)J_n^2(t)dt\right]=
\E\left[\sum_{n=1}^{\sigma-1}\int_{0}^{1}(1-t)J_n^2(t)dt\right]
\defeg C \in \R.
\end{equation*}
Furthermore we have
\begin{eqnarray*}
\left|\E\left[\sum_{n=1}^{\sigma-1}e^{-\lambda \sum_{k=0}^{n} Z_k}
\!\int_{0}^{1}\!(1-t)J_n^2(t)dt\right]-C\right|&\!\!\!=\!\!\!&
\left|\E\left[\sum_{n=1}^{\sigma-1}(1-e^{-\lambda \sum_{k=0}^{n}
Z_k})
\!\int_{0}^{1}\!(1-t)J_n^2(t)dt\right]\right|\\
&\!\!\!\le\!\!\! & \Cste{43} \E\left[(1-e^{-\lambda
\sum_{k=0}^{\sigma-1} Z_k})\sum_{n=1}^{\sigma-1}
Z_n^{\frac{\alpha}{2}}\right].
\end{eqnarray*}
And using Hölder's inequality, we get
\begin{eqnarray*}
 \E\left[(1-e^{-\lambda \sum_{k=0}^{\sigma-1}
Z_k})\sum_{n=1}^{\sigma-1} Z_n^{\frac{\alpha}{2}}\right]&\le &
\E\left[(1-e^{-\lambda \sum_{k=0}^{\sigma-1}
Z_k})\sigma^{\frac{1}{3}}(\sum_{n=1}^{\sigma-1}
Z_n^{\frac{3\alpha}{4}})^{\frac{2}{3}}\right]\\
&\le & \E\left[(1-e^{-\lambda \sum_{k=0}^{\sigma-1}
Z_k})^3\sigma\right]^{\frac{1}{3}}\E\left[\sum_{n=1}^{\sigma-1}
Z_n^{\frac{3\alpha}{4}}\right]^{\frac{2}{3}}\\
&\le &\Cste{44} \E\left[(1-e^{-\lambda \sum_{k=0}^{\sigma-1}
Z_k})\sigma\right]^{\frac{1}{3}}\\
 &\leq& \Cste{45}\lambda^\varepsilon
\end{eqnarray*}
where we used Lemma \ref{ART5_lemmutil} for the last inequality. This
yields (\ref{ART5_J3proof}) and completes, at last, the proof of (e) of
Lemma \ref{ART5_equivmun} when $\alpha\in (0,1)$.

\subsubsection{Proof of  Lemma \ref{ART5_equivmun} when $\alpha=1$}
The proof of the lemma when $\alpha=1$ is quite similar to the one for $\alpha<1$.
Giving a complete proof would be quite lengthy and
redundant. We shall therefore provide only the arguments which
differ from the case $\alpha<1$.

For $\alpha=1$, the main difference from the previous case comes
from the fact that the function $F_{1-\nu}=F_0$ is not bounded near
0 anymore, a property that was extensively used in the course of the
proof when $\alpha<1$. To overcome this new difficulty, we introduce
the function $G$ defined by
\begin{equation}\label{ART5_defG}
G(x)\defeg F_0(x)+F_1(x)\log x\qquad\hbox{for $x>0$}.
\end{equation}
Using the properties of $F_0$ and $F_1$ stated in section
\ref{ART5_sectionchoixphi}, we easily check that the function $G$
satisfies
\begin{enumerate}
\item[(1)] $G(0)\defeg \lim_{x\to 0^+}G(x) = \log(2)-\gamma$ (where $\gamma$ denotes Euler's constant).
\item[(2)] There exists $c_G>0$ such that $G(x)\leq c_G e^{-x}$ for all $x\ge
0$.
\item[(3)] $G'(x)=-xF_0(x)\log x$, so $G'(0)=0.$
\item[(4)] There exists $c_{G'}>0$ such that $|G'(x)|\le c_{G'}
\sqrt{x}e^{-x/2}$ for all $x\ge 0$.
\end{enumerate}
Thus, each time we encounter $F_0(x)$ in the study of $\mu_k(n)$, we
will write $G(x)-F_1(x)\log x$ instead. Let us also notice that
$F_1$ and $F_1'$ are also bounded on $[0,\infty)$.

We now point out, for each assertion (a) - (e) of Lemma
\ref{ART5_equivmun}, the modification required to handle the case
$\alpha=1$.

\bigskip
\noindent\underline{\textbf{Assertion (a)}: $\E[\mu(0)]=C_0\lambda
\log{\lambda}+C_0'\lambda+o(\lambda)$}
\medskip

\noindent As in section \ref{ART5_preuvea}, we have
\begin{eqnarray*}\E[\mu(0)] &=&\lambda\E\left[\int_0^{Z_1}x
F_{0}(\sqrt{\lambda}x)dx\right]\\&=&\lambda\E\left[\int_0^{Z_1}x
G(\sqrt{\lambda}x)dx\right]-\lambda\E\left[\int_0^{Z_1}x
F_{1}(\sqrt{\lambda}x)\log(\sqrt{\lambda}x)dx\right]\\
&=&\lambda\E\left[\int_0^{Z_1}x
\left(G(\sqrt{\lambda}x)-F_{1}(\sqrt{\lambda}x)\log
x\right)dx\right]-\frac{1}{2}\lambda\log
\lambda\E\left[\int_0^{Z_1}x F_{1}(\sqrt{\lambda}x)dx\right]
\end{eqnarray*}
and by dominated convergence,
\begin{equation*}
\lim_{\lambda\rightarrow 0}\E\left[\int_0^{Z_1}x
\left(G(\sqrt{\lambda}x)-F_{1}(\sqrt{\lambda}x)\log
x\right)dx\right]=\E\left[\int_0^{Z_1}x \Big(G(0)-F_{1}(0)\log
x\Big)dx\right].
\end{equation*}
Furthermore, using the fact that $F_1'$ is bounded, we get
\begin{equation*}
\E\left[\int_0^{Z_1}x
F_{1}(\sqrt{\lambda}x)dx\right]=\frac{F_1(0)}{2}\E[Z_1^2]+\mathcal{O}(\sqrt{\lambda})
\end{equation*}
so that
\begin{equation*}
\E[\mu(0)]=C_0\lambda \log{\lambda}+C_0'\lambda+o(\lambda).
\end{equation*}

\bigskip
\noindent\underline{\textbf{Assertion (b)}:
$\E[\sum_{n=1}^{\sigma-1}\mu_1(n)]=o(\lambda)$}
\medskip

\noindent This result is the same as when $\alpha<1$, the only
difference being that now
\begin{equation*} \P\{Z_\infty>x\}\underset{x \rightarrow
\infty}\sim \frac{C\log x}{x}.
\end{equation*}
Thus, equality (\ref{ART5_xi12}) becomes
\begin{equation*}\lambda^2\E\left[\sum_{n=1}^{\sigma-1}
Z_n^2\Indic_{\{Z_n\le\frac{-2\log
\lambda}{\sqrt{\lambda}}\}}\right]\underset{\lambda\to 0^+}{\sim}
\Cste{46} \lambda^{\frac{3}{2}}|\log \lambda|^{2}
\end{equation*}
and the same upper bound holds.

\bigskip
\noindent\underline{\textbf{Assertion (c)}:
$\E[\sum_{n=1}^{\sigma-1}\mu_2(n)]=C_2\lambda
\log{\lambda}+C_2'\lambda+o(\lambda)$}
\medskip

\noindent Using the definition of $G$, we now have
\begin{eqnarray*}\mu_{2}(n)&=&\lambda Z_n F_{0}(\sqrt{\lambda}Z_n)f_1(Z_n)e^{-\lambda
\sum_{k=0}^{n} Z_k}\\
&=&\lambda Z_n f_1(Z_n)e^{-\lambda \sum_{k=0}^{n} Z_k}\left[
\big(G(\sqrt{\lambda}Z_n)-F_1(\sqrt{\lambda}Z_n)\log(Z_n)\big)-\frac{1}{2}
\log \lambda F_{1}(\sqrt{\lambda}Z_n)\right].
\end{eqnarray*}
Since $f_1(x)$ is equal to 0 for $x\ge M-1$,  we get the following
(finite) limit
\begin{multline*}\lim_{\lambda\rightarrow 0}\E\left[
\sum_{n=1}^{\sigma-1}Z_n f_1(Z_n)e^{-\lambda \sum_{k=0}^{n} Z_k}
(G(\sqrt{\lambda}Z_n)-F_1(\sqrt{\lambda}Z_n)\log(Z_n))\right]=\\
\E\left[ \sum_{n=1}^{\sigma-1}Z_n f_1(Z_n)
(G(0)-F_1(0)\log(Z_n))\right].
\end{multline*}
Using the same idea as in (\ref{ART5_inter2}), using also Lemma
\ref{ART5_lemmutil} and the fact that $F'_1$ is bounded, we deduce
that
\begin{equation*}
\E\left[ \sum_{n=1}^{\sigma-1}Z_n f_1(Z_n)e^{-\lambda \sum_{k=0}^{n}
Z_k} F_1(\sqrt{\lambda}Z_n))\right]=\E\left[
\sum_{n=1}^{\sigma-1}Z_n f_1(Z_n)
F_1(0)\right]+o(\lambda^{\varepsilon})
\end{equation*}
which completes the proof of the assertion.

\bigskip \noindent\underline{\textbf{Assertion (d)}:
$\E[\sum_{n=1}^{\sigma-1}\mu_3(n)]=C_3\lambda
\log{\lambda}+C_3'\lambda+o(\lambda)$}
\medskip

\noindent We do not  detail the proof of this assertion since it is
very similar to the proof of (c).

\bigskip
\noindent\underline{\textbf{Assertion (e)}:
$\E[\sum_{n=1}^{\sigma-1}\mu_4(n)]=C_4'\lambda+o(\lambda)$}
\medskip

\noindent It is worth noticing that, when $\alpha=1$, the
contribution of this term is negligible compared to (a) (c) (d) and
does not affect the value of the constant in Proposition
\ref{ART5_progeny}. This differs from the case $\alpha<1$. Recall that
\begin{equation*}
\mu_{4}(n)=-e^{-\lambda \sum_{k=0}^{n} Z_k}\E[\theta_n\;|\;
\mathcal{F}_n],
\end{equation*}
where $\theta_n$ is given by (\ref{ART5_deftheta}). Recall also the
notation $V_{n,t}\defeg Z_n+t(Z_{n+1}-Z_n)$. Just as in
(\ref{ART5_thetaI1I2}), we write
\begin{equation*}
\E[\theta_n\;|\;\mathcal{F}_n]=\int_{0}^{1}(1-t)(I_n^1(t)+I_n^2(t))dt,
\end{equation*}
with
\begin{eqnarray*}
I_n^1(t)&\defeg&\lambda\E\left[(Z_{n+1}-Z_n)^2\Big(F_{1}(\sqrt{\lambda}V_{n,t})-
F_{1}(\sqrt{\lambda}Z_n)\Big)\;\big|\;\mathcal{F}_n \right]\\
I_n^2(t)&\defeg&
-\lambda\E\left[(Z_{n+1}-Z_n)^2(F_{0}(\sqrt{\lambda}V_{n,t})-
 F_{0}(\sqrt{\lambda}Z_n))\;\big|\;\mathcal{F}_n \right].
 \end{eqnarray*}
It is clear that inequality (\ref{ART5_split1}) still holds \emph{i.e.}
\begin{equation*}
|I_n^1(t)| \le\lambda^{\frac{3}{2}}\E\left[|Z_{n+1}-Z_n|^3
\max_{x\in[Z_n,Z_{n+1}]}\sqrt{\lambda}x
F_{0}(\sqrt{\lambda}x)\;\big|\;\mathcal{F}_n \right].
\end{equation*}
In view of the relation
\begin{equation*}
F_{0}(\sqrt{\lambda}x)=G(\sqrt{\lambda}x)-F_1(\sqrt{\lambda}x)\log
x-\frac{1}{2}F_1(\sqrt{\lambda}x)\log \lambda,
\end{equation*}
and with similar techniques to those used in the case $\alpha<1$, we
can prove that
\begin{equation} \Big|\E\left[\sum_{n=1}^{\sigma-1}e^{-\lambda \sum_{k=0}^{n}
Z_k} \int_{0}^{1}(1-t)I_n^1(t)dt\right]\Big|\le \Cste{47}
\lambda^{\frac{9}{8}}|\log \lambda| = o(\lambda).
\end{equation}
It remains to estimate $I_n^2(t)$ which  we now decompose into four
terms:
\begin{equation*}
I_n^2(t)=-\lambda(\tilde{J}_n^1(t)+\tilde{J}_n^2(t)+\tilde{J}_n^3(t)+\tilde{J}_n^4(t)),
\end{equation*}
with
\begin{eqnarray*}
\tilde{J}_n^1(t)&\defeg&\E\left[(Z_{n+1}-Z_n)^2(G(\sqrt{\lambda}V_{n,t})-
 G(\sqrt{\lambda}Z_n))\;|\;\mathcal{F}_n \right]\\
 \tilde{J}_n^2(t)&\defeg&-\frac{1}{2}\log \lambda\E\left[(Z_{n+1}-Z_n)^2(F_1(\sqrt{\lambda}V_{n,t})-
 F_{1}(\sqrt{\lambda}Z_n))\;|\;\mathcal{F}_n \right]\\
\tilde{J}_n^3(t)&\defeg&-\E\left[(Z_{n+1}-Z_n)^2\log Z_n
(F_{1}(\sqrt{\lambda}V_{n,t})-
F_{1}(\sqrt{\lambda}Z_n))\;|\;\mathcal{F}_n \right]\\
\tilde{J}_n^4(t)&\defeg&-\E\left[(Z_{n+1}-Z_n)^2(\log
V_{n,t}-\log(Z_n))F_{1}(\sqrt{\lambda}V_{n,t})\;|\;\mathcal{F}_n
\right].
\end{eqnarray*}
We can obtain an upper bound of order $\lambda^{\varepsilon}$ for
$\tilde{J}_n^1(t)$ by considering again three cases:
\begin{equation*}
\hbox{(1) }\frac{1}{2}Z_n<Z_{n+1}<2Z_n \qquad\hbox{(2) }Z_{n+1}\le
\frac{1}{2}Z_n \qquad\hbox{(3) }Z_{n+1}\ge 2Z_n.
\end{equation*}
For (1), we use that $|G'(x)|\leq c_{G'}\sqrt{x}e^{-x/2}$ for all
$x\ge 0$. We deal with (2) combining Lemma \ref{ART5_largedeviation}
and the fact that $G'$ is bounded. Finally, the case (c) may be
treated by similar methods as those used for dealing with
$J_n^{2}(t)$ in the proof of (e) when $\alpha<1$ (\emph{i.e.} we
separate into two terms according to whether $Z_{n+1}\le
\lambda^{-1/4}$ or not).

Keeping in mind that $F_1$ is bounded and that $|F_1'(x)|=x
F_0(x)\le \Cste{48}\sqrt{x}e^{-x}$,  the same method enables us to
deal with $\tilde{J}_n^2(t)$ and $\tilde{J}_n^3(t)$. Combining these
estimates, we get
\begin{equation*}
\E\left[\sum_{n=1}^{\sigma-1}e^{-\lambda \sum_{k=0}^{n} Z_k}\int_0^1
(1-t)\left(\tilde{J}_n^1(t)+\tilde{J}_n^2(t)+\tilde{J}_n^3(t)\right)dt\right]
= o(\lambda^\varepsilon).
\end{equation*}
for $\varepsilon>0$ small enough. Therefore, it merely remains to
prove that
\begin{equation}\label{ART5_jn4}
\lim_{\lambda\to 0^+}\E\left[\sum_{n=1}^{\sigma-1}e^{-\lambda
\sum_{k=0}^{n} Z_k}\int_0^1 (1-t)\tilde{J}_n^4(t)dt\right]
\end{equation}
exists and is finite. In view of the dominated convergence theorem,
it suffices to prove that
\begin{equation}\label{ART5_jn4abs}
\E\left[\sum_{n=1}^{\sigma-1}\int_0^1
(1-t)\E\Big[(Z_{n+1}-Z_n)^2|\log
V_{n,t}-\log(Z_n)|\;\Big|\;\mathcal{F}_n \Big]dt\right]<\infty.
\end{equation}
We consider separately the cases $Z_{n+1}> Z_n$ and $Z_{n+1}\le
Z_n$. On the one hand, using the inequality $\log(1+x)\le x$, we get
\begin{multline*}
\E\Big[\Indic_{\{Z_{n+1}> Z_n\}}(Z_{n+1}-Z_n)^2|\log
V_{n,t}-\log(Z_n)|\;\Big|\;\mathcal{F}_n \Big]\\
\le \E\Big[\Indic_{\{Z_{n+1}> Z_n\}}(Z_{n+1}-Z_n)^2\log
\Big(1+\frac{t(Z_{n+1}-Z_n)}{Z_n}\Big)\;\Big|\;\mathcal{F}_n
\Big]\le t\sqrt{Z_n}.
\end{multline*}
On the other hand, we find
\begin{multline*}
\E\Big[\Indic_{\{Z_{n+1}\le Z_n\}}(Z_{n+1}-Z_n)^2|\log
V_{n,t}-\log(Z_n)|\;\Big|\;\mathcal{F}_n \Big]\\
\le \E\Big[\Indic_{\{Z_{n+1}\le Z_n\}}(Z_{n+1}-Z_n)^2\log
\Big(1+\frac{t(Z_{n}-Z_{n+1})}{Z_n-t(Z_n-Z_{n+1})}\Big)\;\Big|\;\mathcal{F}_n
\Big]\le \frac{t}{1-t}\sqrt{Z_n}.
\end{multline*}
Since $\E[\sum_{n=1}^{\sigma-1}\sqrt{Z_n}]$ is finite,  we deduce
 (\ref{ART5_jn4abs}) and the proof of assertion (e) is complete.

\section{Proof of Theorem \ref{ART5_MainTheo}}\label{ART5_application}
Recall that $X$ stands for the $(M,\p)$-cookie random walk and $Z$
stands for its associated branching process. We define the sequence
of return times $(\sigma_n)_{n\geq 0}$ by
\begin{equation*}
\left\{
\begin{array}{rll}
\sigma_0 &\defeg& 0,\\
\sigma_{n+1} &\defeg& \inf\{k>\sigma_n\, ,\, Z_k=0\}.
\end{array}
 \right.
\end{equation*}
In particular, $\sigma_1 = \sigma$ with the notation of the previous
sections. We write
\begin{equation*}
\sum_{k=0}^{\sigma_n}Z_k=\sum_{k=\sigma_0}^{\sigma_1-1}Z_k+\ldots+\sum_{k=\sigma_{n-1}}^{\sigma_n-1}Z_k.
\end{equation*}
The random variables $(\sum_{k=\sigma_i}^{\sigma_{i+1}-1}Z_k\, , \,
i\in\N)$ are i.i.d. In view of Proposition \ref{ART5_progeny}, the
characterization of the domains of attraction to a stable law
implies
\begin{equation}\label{ART5_fin1}
\left\{\begin{array}{lll}
\frac{\sum_{k=0}^{\sigma_n}Z_k}{n^{1/\nu}}&\underset{n\to\infty}{\overset{\hbox{\tiny{law}}}{\longrightarrow}}\mathcal{S}_{\nu}
& \hbox{ when } \alpha\in(0,1),\\
\frac{\sum_{k=0}^{\sigma_n}Z_k}{n\log
n}&\underset{n\to\infty}{\overset{\hbox{\tiny{prob}}}{\longrightarrow}}
c & \hbox{ when } \alpha=1.
\end{array}\right.
\end{equation}
where $\mathcal{S}_\nu$ denotes a positive, strictly stable law with
index $\nu \defeg \frac{\alpha+1}{2}$ and where $c$ is a strictly
positive constant. Moreover, the random variables
$(\sigma_{n+1}-\sigma_n\, ,\, n\in\N)$ are i.i.d. with finite
expectation $\E[\sigma]$, thus
\begin{equation}\label{ART5_fin2}
\frac{\sigma_n}{n}\overset{\hbox{\tiny{a.s.}}}{\underset{n\to\infty}{\longrightarrow}}\E[\sigma].
\end{equation}
The combination of (\ref{ART5_fin1}) and (\ref{ART5_fin2}) easily gives
\begin{equation*}
\left\{\begin{array}{lll}
\frac{\sum_{k=0}^{n}Z_k}{n^{1/\nu}}&\underset{n\to\infty}{\overset{\hbox{\tiny{law}}}{\longrightarrow}}\E[\sigma]^{-\frac{1}{\nu}}\mathcal{S}_{\nu} & \hbox{ when } \alpha\in(0,1),\\
\frac{\sum_{k=0}^{n}Z_k}{n\log
n}&\underset{n\to\infty}{\overset{\hbox{\tiny{prob}}}{\longrightarrow}}
c \E[\sigma]^{-1} & \hbox{ when } \alpha=1.
\end{array}\right.
\end{equation*}
Concerning the hitting times of the cookie random walk $T_n =
\inf\{k\geq 0\, ,\, X_k=n\}$, making use of Proposition
\ref{ART5_propZT}, we now deduce that
\begin{equation*}\left\{\begin{array}{lll}
\frac{T_n}{n^{1/\nu}}&\underset{n\to\infty}{\overset{\hbox{\tiny{law}}}{\longrightarrow}}2\E[\sigma]^{-\frac{1}{\nu}}\mathcal{S}_{\nu} & \hbox{ when } \alpha\in(0,1),\\
\frac{T_n}{n\log
n}&\underset{n\to\infty}{\overset{\hbox{\tiny{prob}}}{\longrightarrow}}
2c\E[\sigma]^{-1} & \hbox{ when } \alpha=1.
\end{array}\right.\end{equation*}
Since $T_n$ is the inverse of $\sup_{k\leq n}X_k$, we conclude that
\begin{equation*}
\left\{\begin{array}{lll}
\frac{1}{n^{\nu}}\sup_{k\leq n}X_k&\underset{n\to\infty}{\overset{\hbox{\tiny{law}}}{\longrightarrow}}\mathcal{M}_\nu & \hbox{ when } \alpha\in(0,1),\\
\frac{\log n}{n}\sup_{k\leq
n}X_k&\underset{n\to\infty}{\overset{\hbox{\tiny{prob}}}{\longrightarrow}}
C & \hbox{ when } \alpha=1,
\end{array}\right.\end{equation*}
where $C \defeg (2c)^{-1}\E[\sigma] >0$ and $\mathcal{M}_\nu \defeg
2^{-\nu}\E[\sigma] \mathcal{S}_\nu^{-\nu}$ is a Mittag-Leffler
random variable with index $\nu$. This completes the proof of the
theorem for $\sup_{k\leq n}X_k$. It remains  to prove that this
result also holds for $X_n$ and for $\inf_{k\ge n}X_k$. We need the
following lemma.

\begin{lemm}\label{ART5_supmoinsinf}Let $X$ be a transient cookie random walk. There exists $f:\N\mapsto \R_+$ with
$\lim_{K\rightarrow +\infty}f(K)=0$ such that, for every $n\in \N$,
$$\P\{n-\inf_{i\ge T_n}X_i>
K\}\le f(K).$$
\end{lemm}

\begin{proof}
The proof of this lemma is very similar to that of Lemma 4.1 of
\cite{BasdevantSingh06-preprint}. For $n \in \N$, let $\omega_{X,n} =
(\omega_{X,n}(i,x))_{i\geq 1,x\in\Z}$ denote the random cookie
environment at time $T_n$ "viewed from the particle", \emph{i.e.} the
environment obtained at time $T_n$ and shifted by $n$. With this
notation, $\omega_{X,n}(i,x)$ denotes the strength of the
$i^{\hbox{\tiny{th}}}$ cookies at site $x$:
\begin{equation*}
\omega_{X,n}(i,x)= \left\{ \begin{array}{ll} p_j & \hbox{if }
j=i+\sharp\{0\le k<T_n, X_k=x+n\} \leq M,\\
\frac{1}{2}&\hbox{otherwise.}
\end{array}
\right.
\end{equation*}
Since the cookie random walk $X$ has not visited the half line
$[n,\infty)$ before time $T_n$, the cookie environment
$\omega_{X,n}$ on $[0,\infty)$ is the same as the initial cookie
environment, that is, for $x\geq 0$,
\begin{equation}\label{ART5_fin3}
\omega_{X,n}(i,x)= \left\{
\begin{array}{ll}
p_i&\hbox{if }1\leq i \leq M,\\
\frac{1}{2}&\hbox{otherwise.}
\end{array}
\right.
\end{equation}
Given a cookie environment $\omega$, we denote by $\P_\omega$ a
probability under which $X$ is a cookie random walk starting from
$0$ in the cookie environment $\omega$. Therefore, with these
notations,
\begin{equation}\label{ART5_environmentwn}\P\{n-\inf_{i\ge T_n}X_i> K\}\le \E\left[\P_{\omega_{X,n}}\{X
\mbox{ visits $-K$ at least once}\}\right].\end{equation} Consider
now the deterministic (but non-homogeneous) cookie environment
$\omega_{\p,+}$ obtained from the classical homogeneous $(M,\p)$
environment by removing all the cookies situated on $(-\infty,-1]$:
\begin{equation*} \left\{
\begin{array}{l}
\omega_{\p,+}(i,x)=\frac{1}{2}, \quad \mbox{ for all } x< 0 \mbox{
and }
i\ge 1, \\
\omega_{\p,+}(i,x)=p_i, \quad \mbox{ for all } x\ge 0 \mbox{ and }
i\ge 1 \hbox{ (with the convention $p_i=\frac{1}{2}$ for $i\geq
M$).}
\end{array}\right.
\end{equation*}
According to (\ref{ART5_fin3}), the random cookie environment
$\omega_{X,n}$ is almost surely larger than the environment
$\omega_{\p,+}$ for the canonical partial order, \emph{i.e.}
\begin{equation*}
\omega_{X,n}(i,x)\ge \omega_{\p,+}(i,x)\quad\hbox{for all $i\geq 1$,
$x\in \Z$, almost surely.}
\end{equation*}
The monotonicity result of Zerner stated in Lemma $15$ of
\cite{Zerner05} yields
\begin{equation*}
\P_{\omega_{X,n}}\{X \mbox{ visits } -K\mbox{ at least once}\} \leq
\P_{\omega_{\p,+}}\{X \mbox{ visits } -K\mbox{ at least once}\}
\quad\mbox{almost surely.}
\end{equation*}
Combining this with (\ref{ART5_environmentwn}), we get
\begin{equation}\label{ART5_fin4}
\P\{n-\inf_{i\ge T_n}X_i> K\}\le \P_{\omega_{\p,+}}\{X \mbox{ visits
} -K\mbox{ at least once}\}.
\end{equation}
This upper bound does not depend on $n$. Moreover, it is shown
in the proof of Lemma $4.1$ of \cite{BasdevantSingh06-preprint} that the
walk in the cookie environment $\omega_{\p,+}$ is transient which
implies, in particular,
\begin{equation*}
\P_{\omega_{\p,+}}\{X \mbox{ visits } -K\mbox{ at least once}\}
\underset{K\to\infty}{\longrightarrow}0.
\end{equation*}
\end{proof}

We now complete the proof of  Theorem \ref{ART5_MainTheo}. Let
$n,r,p \in \N$, using the equality $\{T_{r+p}\leq n\}=\{\sup_{k\le
n}X_k\geq r+p\}$, we get
\begin{equation*}
\{\sup_{k\le n}X_k < r\}\subset \{\inf_{k\ge n}X_k < r\}\subset
\{\sup_{k\le n}X_k < r+p\}\cup \{\inf_{k\ge T_{r+p}} X_k< r\}.
\end{equation*}
Taking the probability of these sets, we obtain
\begin{equation*}
\P\{\sup_{k\le n}X_k < r\}\le \P\{\inf_{k\ge n}X_k < r\}\le
\P\{\sup_{k\le n}X_k < r+p\}+ \P\{\inf_{k\ge T_{r+p}} X_k< r\}.
\end{equation*}
But, using Lemma \ref{ART5_supmoinsinf}, we have
\begin{equation*}
\P\{\inf_{k\ge T_{r+p}} X_k< r\}=\P\{r+p-\inf_{k\ge T_{r+p}} X_k>
p\}\le f(p) \underset{p\to\infty}{\longrightarrow} 0.
\end{equation*}
Choosing $x\geq 0$ and $r=\lfloor xn^\nu\rfloor$ and
$p=\lfloor\log n\rfloor$, we get, for $\alpha<1$, as $n$ tends to infinity
$$\lim_{n\rightarrow \infty}\P\left\{\frac{\inf_{k\ge n}X_k}{n^\nu}<
x\right\}=\lim_{n\rightarrow \infty}\P\left\{\frac{\sup_{k\le n}X_i}{n^\nu}< x\right\} = \P\left\{\mathcal{M}_\nu < x\right\}.$$
Of course, the same method also works when $\alpha=1$. This proves
Theorem \ref{ART5_MainTheo} for $\inf_{k\ge n}X_k$. Finally, the result for $X_n$ follows from
$$\inf_{k\ge n}X_k\le X_n \le \sup_{k\le n}X_k.$$

\begin{ack} The authors wish to thank Yueyun Hu whose help, advices and ideas were essential for the development of this paper.
\end{ack}
\bibliographystyle{plain}
\bibliography{Cookie-growth}

\end{document}